\documentclass[]{article}   
\usepackage{amsmath,amsfonts,amsthm,amssymb,dsfont,mathtools}

\usepackage{amssymb}
\usepackage{graphicx}
\usepackage{enumerate}

\usepackage[normalem]{ulem}
\usepackage[dvipsnames,svgnames]{xcolor}
\usepackage{fancyvrb}
\usepackage{subfigure}
\usepackage{epstopdf} 

\usepackage{tikz}

\theoremstyle{plain}
\newtheorem{theorem}{Theorem}[section]

\newtheorem{proposition}[theorem]{Proposition}
\newtheorem{corollary}[theorem]{Corollary}

\theoremstyle{definition}
\newtheorem{definition}[theorem]{Definition}
\newtheorem{example}[theorem]{Example}

\theoremstyle{remark}
\newtheorem{remark}[theorem]{Remark}

\def\CC{\hbox{C\kern -.58em {\raise .54ex \hbox{$\scriptscriptstyle |$}}
  \kern-.55em {\raise .53ex \hbox{$\scriptscriptstyle |$}} }}

\def\eg{{\sl \thinspace e.g.},\ }
\def\ie{{\sl \thinspace i.e.},\ }
\def\pt{\hbox{\bf.}}
\def\ZZ{{{\rm Z}\kern-.28em{\rm Z}}}
\def\RR{\mathop{{\rm I}\kern-.2em{\rm R}}\nolimits}
\def\un{\hbox{\rm1\kern-.28em\hbox{I}}}

\newcommand{\E}{\mathbb{E}}

\newcommand{\F}{\mathbb{F}}
\newcommand{\PP}{\mathbb{P}}

\newcommand{\T}{\mathbb{T}}
\newcommand{\U}{\mathbb{U}}

\renewcommand{\leq}{\leqslant}
\renewcommand{\geq}{\geqslant}
\renewcommand{\phi}{\varphi}
\renewcommand{\Sigma}{\varSigma}
\newcommand{\id}{{\hbox{1\kern-.28em\hbox{I}}}}

\newcommand{\e}{\varepsilon}

\def\un{\hbox{\rm1\kern-.28em\hbox{I}}}

\title{Critical length: an alternative approach
}
\author{
Carolina Vittoria Beccari$^1$, Giulio Casciola$^1$, Marie-Laurence Mazure$^2$}
\date{\it \small
$^1$ Department of Mathematics, University of Bologna, P.zza di Porta San Donato 5, \\
40126 Bologna, Italy\\
$^2$ Universit\'e Grenoble Alpes,  Laboratoire Jean Kuntzmann, 
CNRS, UMR 5224, \\
BP 53, F-38041 Grenoble 9, France \\
\rm carolina.beccari2@unibo.it, giulio.casciola@unibo.it,\\
marie-laurence.mazure@univ-grenoble-alpes.fr
}

\begin{document}
\thispagestyle{empty}
\maketitle
\thispagestyle{empty}
\pagestyle{plain}

\begin{abstract}
We provide a numerical method to determine the critical lengths of  linear differential operators with constant real coefficients.  The need for such a procedure arises when the orders increase. The interest of this article is clearly on the practical side since knowing the critical lengths permits an optimal use of the associated kernels. The efficiency of the procedure is due to its being based on crucial features of Extended Chebyshev spaces on closed bounded intervals. 
\end{abstract}


\section{Introduction}

 The {\it critical length} is a crucial notion attached to kernels of linear differential operators  with constant real coefficients, which was first introduced in \cite{CMP04}. 
Such kernels are known to advantageously replace polynomial spaces in many situations. This is due to the fact that, unlike polynomial spaces,  they inherently depend on parameters which can be used to 
modify the solution to classical problems  
(e.g., interpolation, design, approximation). Nevertheless, to take benefit of these parameters it may be necessary to restrict the length of the interval $[a,b]$ we are working on.  
For instance, for a given differential operator $L$,  if we are interested in Hermite interpolation, $\E_L:=\ker L$ must be an {\it Extended Chebyshev space on $[a,b]$}. This is ensured if and only if the length $b-a$ is less than a fixed number $\ell_L \in ]0, +\infty]$. This number $\ell_L$ is referred to as the {\it critical length} of $L$ (or of $\E_L$).  
If we want to use $\E_L$ for design \cite{CA05}, we have to require $\E_L$ to contain the constants and the length $b-a$ to be less than the critical length of the space $D\E_L$ obtained by differentiation, which is less than or equal to the critical length of $\E_L$ and is called the {\it critical length for design}  of $\E_L$.

It is well known that $\ell_L =+\infty$ if and only if the characteristic polynomial of $L$ has only real roots. Therefore, determining  the critical lengths concerns only  differential operators whose characteristic polynomials have at least one non-real root. The classical approach consists in finding the smallest positive zero of a number of Wronskians attached to $L$ \cite{CMP04, JJA12}.  
Unfortunately,  this Wronskian approach is generally 
difficult to carry out in practice, all the more so as the  dimension increases. As an example, consider the simplest case of cycloidal spaces (i.e., spaces spanned by  polynomials of some degree and the two functions $\cos$ and $\sin$) commonly used in geometric design, see, \eg \cite{HPW, zhang, HPML, beyond, MPS-R, ChenWang,  hoff06, MP}. Their critical lengths were successively  investitaged  in a series of articles \cite{CMP04, CMP14, CMP17}.  They  were definitely identified in relation to  zeros of Bessel functions in \cite{CMP17}.  It is worthwhile mentioning a surprising  feature attached to this class of spaces: the critical length increases only every two steps of dimension elevation. This property is connected with the fact that, whatever the dimension, only one Wronskian is really involved in the computation of the critical length. 

Cycloidal spaces are among the spaces which have also been extensively used during the last decade to build {\it generalised splines} (which are themselves examples of {\it Chebyshevian} or {\it piecewise Chebyshevian splines}) for Isogeometric Analysis purposes \cite{CosManPelSam, ManPelSam, ManPelSpe, ManReaSpe, IGA}. They  belong to the  larger class of spaces which are closed under reflection, corresponding to either odd or even characteristic polynomials. 
Such spaces were fully investigated in  dimension four in \cite{BM12}, with both the exact computation of all the critical lengths, and the analysis of the shape effects they produce in dimension five, within therefore the critical lengths for design. This study was motivated by the fact that dimension four is the lowest dimension in which hyperbolic and trigonometric functions can be combined within the same space. 

Obviously, the critical length should be known prior to any use of a space requiring to work with Extended Chebyshev spaces. To the best of our knowledge, no other computation of critical lengths exist apart from the examples we mentioned, except for  trivial cases (\eg lower dimensions, or spaces resulting from changes of variables). Beyond dimension four an exact computation is not really expectable anyway. This encouraged us to develop an effective numerical procedure instead. An advantage of the proposed algorithm is that it simultaneously provides not only the Bernstein-type bases which can then be used for numerical computations, but also associated generalised derivatives which can serve to develop approximation properties \cite{schu}. 

Our numerical procedure is described and illustrated in Section 4. It is entirely based on the question: {\it how to connect a number of Extended Chebyshev spaces on consecutive intervals, all of the same dimension, so as to produce a global Extended Chebyshev space of the same dimension?} A numerical answer to this question is briefly presented in Section 3, as an application  of a more general numerical test including connection matrices which was developed in \cite{BCM17}. Beforehand, Section 2 provides readers with the necessary background on which the test  relies. We more specifically underline the importance of dimension diminishing through generalised derivatives, and its action  on Bernstein and Bernstein-like bases, which plays a prominent role  in the present work.   Section 5 concludes with a few crucial comments drawing readers' attention on why it is both necessary and advantageous to know the critical length. Advantageous: only the knowledge of its critical length (for design) enables  us  to take the maximum benefit of a given space, in particular in view of  constructing splines with pieces taken from  different spaces. Necessary: in terms of design, for instance, visual shape preservation is certainly not sufficient,  the only safe approach being to work within critical length for design. Not only do these observations show the interest of our numerical approach,  but they are all the more crucial as they are somehow in contradiction with rather widely spread practices. Last but not least, we would like to mention that our numerical procedure can also provide  useful help to raise conjectures / solve theoretical questions. To cite only one point, it  clearly indicates that  the every two step increase, or the involvement of only one Wronskian, observed for cycloidal spaces, concerns only a limited subclass of spaces.

\section{Background}

In this section we gather the basic concepts and properties concerning Extended Chebyshev spaces strongly involved in the present paper. For further acquaintance with the subject, see \cite{karlin, schu, HP, CA99, ready} for instance.
\subsection{Extended Chebyshev spaces}

Let  $I$ be a non-trivial real interval and let $\E_n\subset C^n(I)$ be an $(n+1)$-dimensional space. Then, $\E_n$ is said to be an {\it Extended Chebyshev space on $I$} (for short, EC-space on $I$), if any non-zero $F\in \E_n$ vanishes at most $n$ times in $I$, counting multiplicities up to $(n+1)$, or equivalently, if any Hermite interpolation problem in $(n+1)$ data in $I$ is unisolvent in $\E_n$. It is said to be a {\it W-space on $I$}, if the Wronskian of any basis of $E_n$ never vanishes in $I$, or equivalently, if any Taylor interpolation problem in $(n+1)$ data in $I$ is unisolvent in $\E_n$. An $(n+1)$-dimensional EC-space on $I$ is thus a W-space on $I$ but the converse property is not true, except for $n=0$. 

As is well known, the class of all W-spaces on $I$ is closed under integration and multiplication by a sufficiently differentiable function which does not vanish on $I$, and the same holds true for the class of all EC-spaces on $I$. 

While it is inherent in their definition that EC-spaces are crucial for interpolation, in Theorem \ref{ECfordesign} below we remind the reader why they are crucial for design too. Beforehand let us recall below the definition of bases of the Bernstein-type. 
\begin{definition}
Given $a,b\in I$, $a<b$, we say that a sequence $(V_0, \ldots, V_n)$ of functions in $C^n(I)$ is  {\it a Bernstein-like basis relative to $(a,b)$} if, for $i=0, \ldots, n$, the function $V_i $ vanishes exactly $i$ times at $a$ and exactly $(n-i)$ times at $b$.  A {\it positive Bernstein-like basis relative to $(a,b)$} is a Bernstein-like basis $(V_0, \ldots, V_n)$ relative to $(a,b)$ such that $V_i$ is positive on $]a,b[$ for $i=0, \ldots, n$.
\end{definition}

\begin{definition}
Given $a,b\in I$, $a<b$, a Bernstein basis relative to $(a,b)$ is a positive Bernstein-like basis $(B_0, \ldots, B_n)$ relative to $(a,b)$ which is normalised, \ie $\sum_{i=0}^nB_i=\un$, where $\un$ stands for the constant function $\un(x)=1$ for all $x\in I$. 
\end{definition}

Most of the time, for the sake of simplicity, it is convenient to include the positivity in the terminology ``Bernstein-like basis", but in the present paper it is essential to state it separately. Indeed, proving the positivity of the various bases is a major concern in the numerical test described in the next section. From now on, $D$ denotes the ordinary differentiation on any interval.

\begin{theorem}
\label{ECfordesign} 
For a given $(n+1)$-dimensional W-space $\E_n\subset C^n(I)$, supposed to contain the constants, the following properties are equivalent:
\begin{enumerate}[\rm(\roman{enumi})]
\item for each $a,b\in I$, $a<b$, $\E_n$ possesses a Bernstein basis relative to $(a,b)$;
\item for each $a,b\in I$, $a<b$, $D\E_n$ possesses a Bernstein-like basis relative to $(a,b)$;
\item the ($n$-dimensional) space $D\E_n$ is an EC-space on $I$;
\item blossoms exist in the space $\E_n$. 
\end{enumerate}
Furthermore, when {\rm(iii)} is satisfied, all the classical design algorithms can be developed in $\E_n$, and for each $a,b\in I$, $a<b$, the Bernstein basis relative to $(a,b)$ is the optimal normalised totally positive basis in $\E_n$ restricted to $[a,b]$. 
\end{theorem}

We mention blossoms only because they are the underlying tool for many of the results involved in the present work. We will not say more on them. Readers interested can refer to \cite{HP} and to many articles by the third author. Assuming that (ii) holds true, let $(B_0, \ldots, B_n)$ be the Bernstein basis in $\E_n$, relative to $(a,b)\in I^2$, $a<b$. Its total positivity on $[a,b]$ means that, for any $a\leq x_0<x_1<\cdots <x_n\leq b$, all minors of the matrix  $\bigl(B_i(x_j)\bigr)_{0\leq i,j\leq n}$ are non-negative. This is known to guarantee shape preserving properties in $\E_n$, see \cite{TNT}. The mentioned optimality  refers to the fact that we cannot find a better basis regarding this question, see \cite{optimal, CMP04}.  These comments justify the following definition:
\begin{definition} 
\label{def:GFD}
For any $n\geq 1$, an $(n+1)$-dimensional W-space $\E_n$ on $I$ is said to be {\it good for design} when
first,  it contains the constants, and second, the space $D\E_n$ is an EC-space on $I$. 
\end{definition}

\noindent
Observe that a  W-space $\E_n$ which  is good for design is automatically an EC-space on $I$. 

\begin{remark}The closure of the class of all EC-spaces on $I$ under multiplication by positive functions and integration can be visualised as follows: 
\begin{equation}
\label{ECfromw}
	\begin{array}{lc}
		\hbox{\bf step 2:}\quad&\E_{k+1}=(k+2)\hbox{-dimensional EC-space on }I\\
	&\hbox{\it multiply by any }\uparrow \hbox{\it positive } w\in C^{k+1}(I)\\
	\hbox{\bf step 1:}\quad&\F_{k+1}=(k+2)\hbox{-dimensional EC-space good for design on }I,
	\\
	&\quad\, \hbox{\it inte-}\uparrow \hbox{\it grate } \\
	\hbox{\bf step 0:}\quad&\E_k=(k+1)\hbox{-dimensional EC-space on }I
\end{array}
\end{equation}
Select a sequence $(w_0, \ldots, w_n)$ of {\it weight functions on $I$},  in the sense that, for $i=0, \ldots, n$ $w_i$ is $C^{n-i}$ and positive on $I$.  For $k=0, \ldots, n$, we can then repeatedly apply the dimension increasing scheme (\ref{ECfromw}) starting from the space $\E_0$ spanned by $w_n$, where the passage from $\F_k$ to $\E_k$ corresponds to multiplication by $w_{n-k}$, $k=0, \ldots, n$.  The whole process takes place within the class of all EC-spaces on $I$. If we denote  by $EC(w_0, \ldots, w_n)$ the final EC-space  $\E_n$, the same notation yields
$$
\E_k=EC(w_{n-k}, \ldots, w_n), \quad \F_k=EC(\un, w_{n-k+1}, \ldots, w_n),\  k=0, \ldots, n, 
$$
in which we have added the space $\F_0$ of all constant functions on $I$, and the space $\F_{n+1}$ obtained by integration of $\E_n$. Note that each space $\F_k$, $k=1, \ldots, n+1$, is good for design on $I$. Classically, the system $(w_0, \ldots, w_n)$ is associated with linear differential operators $L_0, \ldots, L_n$ -- also named {\it generalised derivatives} -- recursively defined as follows: 
\begin{equation}
\label{difop}
L_0F:=\frac{F}{w_0}\ ,\quad L_iF:=\frac{1}{w_i}DL_{i-1}F,\quad 1\leq i\leq n.
\end{equation}
With these notations, $\E_n=EC(w_0, \ldots, w_n)$ can be described as the set of all $F\in C^n(I)$ such that $L_nF\in\F_0$. 

\medskip

The previous process relates a classical way to obtain EC-spaces on a given interval, the generalised derivatives enabling the development of important approximation properties modelled on polynomial spaces \cite{schu}.  In particular, not only does this provide one final EC-space but even a nested sequence of EC-spaces
$$
\E_0\subset \E_1\subset \cdots \subset \E_{n-1}\subset \E_n, \quad \hbox{with }\E_i:=EC(w_0, \ldots, w_i) \hbox{ for }i=0, \ldots, n. 
$$ 
This is the reason why $\E_n=EC(w_0, \ldots, w_n)$ is called the {\it Extended Complete Chebyshev space associated with $(w_0, \ldots, w_n)$}. The presence of such a nested sequence in $\E_n$ is crucial, for instance, to define associated Chebyshevian divided differences and Newton-type expansions for  the solution to any Hermite interpolation problem in $\E_n$. Conversely, if we start with a given nested sequence 
$\E_0\subset \cdots \E_i\subset \cdots \subset \E_n$, where,  for $i=0, \ldots, n$,  $\E_i$ is an $(i+1)$-dimensional W-space on $I$, then it is well known that it is a nested sequence of EC-spaces on $I$.  More precisely, selecting any sequence $(U_0, \ldots, U_n)$ such that $U_i\in \E_i\setminus\E_{i-1}$ for $i=0, \ldots, n$, with $\E_{-1}:=\{0\}$, we have \cite{karlin, ready}
\begin{equation}
\label{wi}
\E_i=EC(w_0, \ldots, w_i), \quad w_i:=\e_i\ \frac{W(U_0,...,U_{i-2})\ W(U_0,..., U_i)}{W(U_0,...,U_{i-1})^2}\ ,\quad 0\leq i\leq n. 
\end{equation}
In (\ref{wi}), $W(U_0,..., U_i)$ denotes the Wronskian of the sequence $(U_0, \ldots, U_i)$ with the convention that $W(\emptyset)=\un$, and $\e_i=\pm$ is chosen so as to ensure the positivity of $w_i$. 

\end{remark}

\subsection{EC-spaces and dimension diminishing}

Is it possible to reverse the process (\ref{ECfromw}) so as to find a system $(w_0, \ldots, w_n)$ of weight functions on $I$ such that $\E_n=EC(w_0, \ldots, w_n)$?  The first step should thus consist in a dimension diminishing procedure as follows: 

\begin{equation}
\label{reverse}
	\begin{array}{lc}
	\hbox{\bf step 0:}\quad&\E_n\subset C^n(I),\ (n+1)\hbox{-dimensional}\\
	&\hbox{\it division by a } \downarrow \hbox{\it  positive } w_0\in\E_n\\
	\hbox{\bf step 1:}\quad&\F_n=L_0\E_n\ (n+1)\hbox{-dimensional}  \hbox{ containing }\un \\
	&\hbox{\it differen-} \downarrow \hbox{\it  tiation}\\
	\hbox{\bf step 2:}\quad&\E_{n-1}=DL_0\E_n\subset C^{n-1}(I)\ n\hbox{-dimensional} \\
\end{array}
\end{equation}
We should therefore first be able to find a positive function in a given EC-space $\E_n$. Unfortunately, this is not always possible, as proved by the famous counterexample of the space $\E_1$ spanned on $I=[0, \pi[$ by the two functions $\cos, \sin$, which is an EC-space on $I$ but which does not contain any non-vanishing function. No equality of the form $\E_1=EC(w_0, w_1)$ is thus expectable. Nonetheless, such an equality is possible by restriction to any $[0, b]$, with $0<b<\pi$, as reminded below. 

\begin{theorem}
\label{th:EC(w)}
Let $[a,b]$, $a<b$, be a closed bounded interval, and let $\E_n\subset C^n([a,b])$ be $(n+1)$-dimensional. Then, the following properties are equivalent:
\begin{enumerate}[\rm(\roman{enumi})]
\item $\E_n$ is an EC-space (resp., an EC-space good for design) on $[a,b]$; 
\item there exists a system $(w_0, \ldots, w_n)$  (resp., $(w_1, \ldots, w_n)$) of weight functions on $[a,b]$ such that 
$\E_n=EC(w_0, \ldots, w_n)$ (resp., $\E_n=EC(\un, w_1, \ldots, w_n)$).
\end{enumerate}
\end{theorem}
  
Given an $(n+1)$-dimensional EC-space $\E_n$ on $[a,b]$, and any positive function $w_0\in\E_n$, we are certain that the $n$-dimensional space $DL_0\E_n$  obtained according to (\ref{reverse}) is a W-space on $[a,b]$, but there is no guarantee that it is an EC-space on $[a,b]$. This is clear from the classical example where $\E_2$ is the three-dimensional EC-space on $[0, 3\pi/2]$ spanned by the functions $\un, \cos, \sin$, for which the space $D\E_2$, spanned by $\cos, \sin$ is an EC-space on $[0, \pi/2[$ but not on $[0, \pi/2]$. As a matter of fact, one step of dimension diminishing within the class of all EC-spaces on $[a,b]$ can only  be done according to the rule specified below  \cite{JAT11}. 

\begin{theorem}
\label{th:w0}
Let $\E_n$ be an $(n+1)$-dimensional EC-space on $[a,b]$, and let $(V_0, \ldots, V_n)$ denote a  positive Bernstein-like  basis relative to $(a,b)$ in $\E_n$. Given a function $w_0=\sum_{i=0}^n \alpha_i V_i\in \E_n$, the following properties are equivalent:
\begin{enumerate}[\rm(\roman{enumi})]
\item $\alpha_0, \ldots,\alpha_n$ are all positive;
\item $w_0$ is positive on $[a,b]$, and if $L_0$ stands for the division by $w_0$, the $(n+1)$-dimensional space $L_0\E$  is an EC-space good for design on $[a,b]$, \ie the $n$-dimensional space $DL_0\E$ is an EC-space on $[a,b]$.
\end{enumerate}
\end{theorem}

It is worthwhile mentioning the following straightforward but crucial consequence of Theorem \ref{th:w0}, to be compared with Theorem \ref{ECfordesign}, see \cite{JAT11}.  

\begin{corollary}
\label{cor:w0}
For a given $(n+1)$-dimensional space $\E_n\subset C^n([a,b])$, $n\geq 1$, known to be an EC-space on $[a,b]$,  the following properties are equivalent:
 \begin{enumerate}
 \item $\E_n$ possesses a Bernstein basis relative to $(a,b)$;
 \item $\E_n$ is good for design on $[a,b]$. 
 \end{enumerate} 
\end{corollary} 
\begin{remark}
\label{effectonbases}
The effect of generalised differentiation on the bases is well known and we recall it here. Assume that (i) of Theorem \ref{th:w0} is satisfied, \par
\smallskip
\noindent
$\bullet$ \underbar{From $\E_n$ to $L_0\E_n$}: Division by $w_0$ yields 
\begin{equation}
\label{BfromV}
\un=\sum_{i=0}^n B_i, \quad\hbox{with }B_i:=\frac{\alpha_iV_i}{w_0} \hbox{ for  }i=0, \ldots, n. 
\end{equation}
Clearly, $(B_0, \ldots, B_n)$ is the Bernstein basis relative to $(a,b)$ in the space $L_0\E_n$ which is an EC-space good for design on $[a,b]$. 
Observe that the space $L_0\E_n$ is completely determined by the equivalence class of the  sequence $(\alpha_0, \ldots, \alpha_n)$ under proportionality. Accordingly, we can build infinitely many different such spaces $L_0\E_n$. 

\smallskip
\noindent
$\bullet$ \underbar{From $L_0\E_n$ to $DL_0\E_n$}: 
In close relation with the Bernstein basis $(B_0, \ldots, B_n)$, it is convenient to introduce the functions 
\begin{equation}
\label{Bstar}
B_i^\star:=\sum_{k=i}^nB_k=\un -\sum_{k=0}^{i-1}B_k, \quad i=0, \ldots, n. 
\end{equation}
For each $i=1, \ldots, n$, the function $B_i^\star$ is characterised by the fact that
$$
B_i^\star\hbox{ vanishes exactly }i \hbox{ times at  }a,  
\hbox{ and  }\un-B_i^\star\hbox{  vanishes exactly }(n-i+1) \hbox{ times at }b. 
$$
For this reason, these functions are named {\it transition functions} in the space $L_0\E_n$, see \cite{NUMA19}. Let us set 
\begin{equation}
\label{V=DBstar}
\overline V_i:=DB_{i+1}^\star=\sum_{k=i+1}^nDB_k=-\sum_{k=0}^iDB_k, \quad i=0, \ldots, n-1.
\end{equation}
Clearly,  $(\overline V_0, \ldots,\overline V_{n-1})$ is a Bernstein-like basis relative to $(a,b)$ in the space $DL_0\E_n$ \cite{CA05}. Moreover, expansions in that basis can easily be derived from expansions in the Bernstein basis in $L_0\E_n$, see \cite{CA05}. From (\ref{V=DBstar}) we can also see that each $\overline V_i$ is positive close to $a$. Accordingly, because the space $DL_0\E_n$ is known to be an EC-space on $[a,b]$ (Theorem \ref{th:w0}), we can conclude that  $(\overline V_0, \ldots, \overline V_{n-1})$ is a positive Bernstein-like basis relative to $(a,b)$ which can be used to iterate the process to construct $w_1$  via Theorem \ref{th:w0},  and so forth up to  one   sequence $(w_0, \ldots, w_n)$ of weight functions such that $\E_n=EC(w_0, \ldots, w_n)$. 

 \end{remark}
 
 \begin{remark}
 Starting again with the $(n+1)$-dimensional EC-space $\E_n$ on $[a,b]$, consider one sequence of spaces obtained by iteration of Theorem \ref{th:w0}: 
 $$
  \E_n^{\{0\}}:=\E_n, \qquad  \E_n^{\{p\}}:=DL_{p-1}\E_n=EC( w_{p}, \ldots, w_n)\hbox{ for }p=1, \ldots, n, 
   $$
 corresponding to one given equality $\E_n=EC(w_0, \dots, w_n)$.
 For each $p=0, \ldots, n$, the space $L_p\E_n=EC(\un, w_{p+1}, \ldots, w_n)$ is an $(n-p+1)$-dimensional EC-space good for design on $[a,b]$ in which we denote by $(B_0^{\{p\}}, \ldots, B_{n-p}^{\{p\}})$ the Bernstein basis relative to $(a,b)$.  One step of dimension diminishing transforms each Bernstein basis into the next one via (\ref{V=DBstar}) and (\ref{BfromV}). These relations can be read  in the reverse way, which yields, for each $p=n, n-1, \ldots, 1$, 	
\begin{equation}
\label{IRR}
	\begin{split}
&B_0^{\{p-1\}}(x)=1-\displaystyle\frac{\int_a^xw_p(t)B_0^{\{p\}}(t)\ dt}{ \int_a^bw_p(t)B_0^{\{p\}}(t)\ dt},\\
	\\
&B_i^{\{p-1\}}(x)=\displaystyle\frac{\int_a^x w_p(t)B_{i-1}^{\{p\}}(t)\ dt}{\int_a^bw_p(t)B_{i-1}^{\{p\}}(t)\ dt}
	-
	\displaystyle\frac{\int_a^x w_p(t)B_i^{\{p\}}(t)\ dt}{\int_a^bw_p(t)B_i^{\{p\}}(t)\ dt}, \quad
	1\leq i\leq n-p,\\
	\\
&B_{n-p+1}^{\{p-1\}}(x)=\displaystyle\frac{\int_a^x w_p(t)B_{n-p}^{\{p\}}(t)\ dt}{\int_a^bw_p(t)B_{n-p}^{\{p\}}(t)\ dt}, 	
	\end{split}
\end{equation}
starting from $B_0^{\{n\}}=\un$. These relations, first obtained in \cite{JJA09} through blossoms, are the analogue of the classical integral recurrence relations for polynomial Bernstein bases. From our comments in Remark \ref{effectonbases}, we know that there are infinitely many essentially different possibilities to go from 
$B_0^{\{n\}}=\un$ up to the Bernstein basis  in $L_0\E_n$ according to (\ref{IRR}). It should be observed that, apart from trivial exceptions, none of them can be considered a practical way to calculate the Bernstein basis in $L_0\E_n$, since to the contrary, they are derived from the latter basis. 
 \end{remark}

\subsection{Global versus local }
\label{GlobalLocal}
 
We start again with an $(n+1)$-dimensional EC-space $\E_n$ on $[a,b]$. In $\E_n$, we consider two positive Bernstein-like bases: the first one relative to $(a,b)$, say   $(V_0, \ldots, V_n)$; the second relative  to $(a^*, b^*)$, say $(V_0^*, \ldots, V_n^*)$, where $a\leq a^*<b^*\leq b$. For short, we refer to them as global positive Bernstein-like basis / local  positive Bernstein-like basis. 

Let us expand the global basis in the local one  as follows
\begin{equation}
V_i=\sum_{r=0}^n\gamma_{i, r}V_r^*, \quad i=0, \ldots, n. 
\end{equation}
All coefficients of these expansions are known to be positive (see \cite{BCM17} and other references therein), except possibly in accordance with the zero conditions of $(V_0, \ldots, V_n)$ at the endpoints, that is, $\gamma_{i,r}=0$ for $0\leq r\leq i-1$ if $a^*=a$, and $\gamma_{i,r}=0$ for $i+1\leq r\leq n$ if $b^*=b$. We refer to this fact as {\it the positi\-vity property of local expansions} of the global basis. 

\smallskip
 Let $\E_n^*$ denote the restriction of $\E_n$ to the interval $[a^*, b^*]$. Take any system $(w_0, \ldots, w_n)$ of weight functions on $[a,b]$ such that $\E_n=EC(w_0, \ldots, w_n)$. By restriction to $[a^*, b^*]$, it  generates a system $(w_0^*, \ldots, w_n^*)$ of weight functions on $[a^*,b^*]$ such that $\E_n^*=EC(w_0^*, \ldots, w_n^*)$.  This means that the successive steps of the corresponding dimension diminishing (\ref{reverse}) can be applied simultaneously in $\E_n$ and in $\E_n^*$, with the two associated sequences of positive Bernstein-like bases $(V_0^{\{p\}},  \ldots, V_{n-p}^{\{p\}})$, and $(V_0^{*\{p\}},  \ldots, V_{n-p}^{*\{p\}})$, respectively deduced from the initial bases $(V_0^{\{0\}},  \ldots, V_{n}^{\{0\}}):=(V_0,  \ldots, V_{n})$, and $(V_0^{*\{0\}},  \ldots, V_{n}^{*\{0\}}):=(V_0^*,  \ldots, V_{n}^*)$. At each step we can expand the global basis $(V_0^{\{p\}},  \ldots, V_{n-p}^{\{p\}})$ in the local basis $(V_0^{*\{p\}},  \ldots, V_{n-p}^{*\{p\}})$ as 
\begin{equation}
V_i^{\{p\}}=\sum_{r=0}^{n-p}\gamma_{i, r}^{\{p\}}V_r^{*\{p\}}, \quad i=0, \ldots, n-p. 
\end{equation}
In this process, for each $p\leq n-1$, the coefficients at level $(p+1)$ can be computed from those of level $p$.  In the special case where the weight functions are taken as
$$
w_p:=V_0^{\{p\}}+ \cdots +V_{n-p}^{\{p\}}, \quad p=0, \ldots, n, 
$$
then we have \cite{BCM17}:
\begin{equation}
\label{local_iterate}
\gamma_{i, r}^{\{p+1\}}= \frac{ \sum_{j=i+1 }^{n-p }\gamma_{j, r+1}^{\{p\}} }{  \sum_{j=0 }^{n-p }\gamma_{j, r+1}^{\{p\}} }
- \frac{ \sum_{j=i+1 }^{n-p }\gamma_{j, r}^{\{p\}} }{  \sum_{j=0 }^{n-p }\gamma_{j, r}^{\{p\}} }, \quad 0\leq i,r\leq n-p-1.
\end{equation}

\section{Building a global  EC space from local EC spaces}

Given an $(n+1)$-dimensional space $\E_n\subset C^n(I)$, how to determine whether or not  $\E_n$ is  an EC-space on $I$? From (ii) of Theorem \ref{ECfordesign}, we know that this consists in checking whether all  determinants
$$
\det\left(\U(x), \ldots, \U^{(i-1)}(x),  \U(y), \ldots, \U^{(j-1)}(y)\right), \quad i, j\geq 0, \ i+j=n+1,
$$ 
never vanish for $x,y\in I$, $x<y$, where  $\U:=(U_0, \ldots, U_n)^T$  and $(U_0, \ldots, U_n)$ is any basis in $\E_n$. Moreover, in case the space $\E_n$ is known to be a W-space on $I$, we only have to consider positive integers $i,j$. Except for small values of $n$, it is not easy to check this by hand  and it is not easy either in the general case to do it numerically.  

From now on we consider a closed bounded interval $[a,b]$, $a<b$. From Theorem \ref{th:EC(w)}, we know that $\E_n$ is an EC-space on $[a,b]$ if and only if, 

	\begin{enumerate}[--]
	\item we can find a positive function $w_0\in \E_n$;
	\item we can find a positive function $w_1\in DL_0 \E_n$;
	\item \dots\quad  \ldots
	\item we can find a positive function $w_n\in DL_{n-1}\E_n$,
	\end{enumerate}
where the notations are according to (\ref{difop}). Now, replacing each ``we can find" by ``can we find?", we have at our disposal an easy  theoretical test: if at some stage the answer to the ``can we find?" question is negative, then the initial space $\E_n$ is not an EC-space on $[a,b]$; if all answers are affirmative,  then $\E_n$ is indeed an EC-space on $[a,b]$ and we can even say that $(w_0, \ldots, w_n)$ is a system of weight functions associated with $\E_n$. Observe that, in case  $\E_n$ is known to be a W-space on $[a,b]$, we only need affirmative answers until  the last but one question.

Clearly, in general this is not a realistic test, since exhibiting such positive functions is more or less like pulling a rabbit out of a hat. Nonetheless, it becomes realistic in the situation  addressed subsequently, where   the positivity of functions is checked through the positivity of the coefficients of some expansions  in appropriate bases. 

\smallskip
Throughout the rest of the present subsection, the interval $[a,b]$, $a<b$, is given along with a  sequence $\T=(t_1, \ldots, t_q)$ of $q\geq 1$ knots interior to $[a,b]$, with
$$
t_0:=a<t_1<\dots<t_q<t_{q+1}:=b,
$$
and a given positive integer $n$. From now on, we change the notations, the index of a space  no longer being related to its dimension, but to its interval. For $k=0, \ldots, q$, $\E_k$ is an $(n+1)$-dimensional EC-space on $[t_k, t_{k+1}]$, in which we select a positive Bernstein-like basis $(V_{k,0}, \ldots, V_{k,n})$ relative to $(t_k, t_{k+1})$, which we refer to as the $k$th local positive Bernstein-like basis. Consider the space $\E$ defined by the two conditions:
\begin{enumerate}[---]
\item $\E\subset C^n([a,b])$;
\item for each $k=0, \ldots, q$, the restriction of $\E$ to $[t_k, t_{k+1}]$ is $\E_k$. 
\end{enumerate}
Clearly, the space $\E$ is an $(n+1)$-dimensional W-space on $[a,b]$, but in general it  is not an EC-space on $[a,b]$.  A numerical answer to  the question ``Is $\E$ an EC-space on $[a,b]$" can be obtained thanks to the  test built in \cite{BCM17} in the larger framework of adjacent  EC-spaces tied by connection matrices. Here, it answers the question {\it ``Starting from $\E$, can we iteratively diminish the dimension as in (\ref{reverse})?"} Subsequently, we briefly recall the main ideas / steps, modelled on the iterative dimension diminishing step described in Theorem \ref{th:w0}.  

\medskip
\noindent
{$\bullet$ \bf Numerical test -- Step 0:}
\medskip

In this step 0, we first answer the question: {\it does the W-space $\E$ possess a Bernstein-like basis relative to $(a,b)$?} With the notations already used above, we have to test if
\begin{equation}
\label{BLB?}
\det\left(\U(a), \ldots, \U^{(i-1)}(a),  \U(b), \ldots, \U^{(j-1)}(b)\right)\not=0\quad\hbox{ for }i, j\geq 1, \ i+j=n+1.
\end{equation}
If the answer to (\ref{BLB?}) is negative, we can state that $\E$ is not an EC-space on $[a,b]$. Supposing that it  is affirmative, let $(V_0^{\{0\}}, \ldots, V_n^{\{0\}})$ be a Bernstein-like basis of $\E$ relative to $(a,b)$, determined by fixing the first non-zero derivatives at $a$ or $b$ to be in accordance with its possible positivity. 

For each $k=0, \ldots, q$, consider the $k$th local expansion of the previous basis, in the sense of the expansion of its restriction to $[t_k, t_{k+1}]$ in the local positive Bernstein-like basis $(V_{k,0}^{\{0\}}, \ldots, V_{k,n}^{\{0\}}):=(V_{k,0}, \ldots, V_{k,n})$: 
\begin{equation}
\label{pos?}
{V_i^{\{0\}} }_{\big\vert  [t_k, t_{k+1}]}=\sum_{r=0}^n\gamma_{i,k,r}^{\{0\}} V_{k,r}^{\{0\}},  \quad i=0, \ldots, n, \ k=0, \ldots, q.
\end{equation}
Answer the question: {\it do all these local expansions satisfy the positivity property reminded in Subsection \ref{GlobalLocal}?} 

If the answer is affirmative, proceed to Step 1. If not, we know that $\E$ is not an EC-space on $[a,b]$, and the test stops.

\medskip
\noindent
{$\bullet$ \bf Numerical test -- Step 1:}
\medskip

Due to the positivity property of all local expansions (\ref{pos?}),  $(V_0^{\{0\}}, \ldots, V_n^{\{0\}})$ is a positive Bernstein-like basis relative to $(a,b)$. 
Take 
\begin{equation}
\label{w0w0k}
w_0:=V_0^{\{0\}} + \cdots + V_n^{\{0\}}, \qquad w_{0,k}:={w_0}_{\big\vert [t_k, t_{k+1}]}, \quad k=0, \ldots, q.
\end{equation}
Then, $w_0$ is positive on $[a,b]$ and each $w_{0,k}$ has positive coordinates in the local positive Bernstein-like basis of $\E_k$. Accordingly, we can simultaneously diminish the dimension 

\begin{enumerate}[---]
\item globally, via $w_0$, within the class of all W-spaces on $[a,b]$, replacing $\E$  by $DL_0\E$; 
\item for each $k=0, \ldots, q$,  locally, via $w_{0,k}$ within the class of all EC-spaces on $[t_k,t_{k+1}]$, replacing $\E_k$  by $DL_0^k\E_k$, where $L_0^k$ is the division by $w_{0,k}$.
\end{enumerate}
Meanwhile, each local positive Bernstein-like basis  $(V_{k,0}^{\{0\}}, \ldots, V_{k,n}^{\{0\}})$ is transformed into a positive Bernstein-like basis $(V_{k,0}^{\{1\}}, \ldots, V_{k,n-1}^{\{1\}})$ of $DL_0^k\E_k$, relative to $(t_k, t_{k+1})$, and the global positive Bernstein-like basis $(V_0^{\{0\}},  \ldots, V_{n}^{\{0\}})$ into a global Bernstein-like basis $(V_0^{\{1\}},  \ldots, V_{n-1}^{\{1\}})$ of $DL_0\E$, relative to $(a,b)$. These transformations follow the procedure described in Remark \ref{effectonbases}. 

The next step consists in checking if the coefficients of all  local expansions of the new global basis satisfy the positivity property, in which case we can continue diminishing the dimension both globally and locally  through 
$$
w_1:=V_0^{\{1\}} + \cdots + V_{n-1}^{\{1\}}, \qquad w_{1,k}:={w_1}_{\big\vert [t_k, t_{k+1}]}, \quad k=0, \ldots, q.
$$
All in all, starting from the positive 
$\gamma_{i,k,r}^{\{0\}}$,  $i=0, \ldots, n, \ k=0, \ldots, q$, Step 1  can be translated into a simple iterative computation of real numbers at level $(p+1)$, $\gamma_{i, k, r}^{\{p+1\}}$, from positive $\gamma_{i, k, r}^{\{p\}}$ at level $p$,  according to the formul\ae:
\begin{equation}
\gamma_{i, k, r}^{\{p+1\}}= \frac{ \sum_{j=i+1 }^{n-p }\gamma_{j, k, r+1}^{\{p\}} }{  \sum_{j=0 }^{n-p }\gamma_{j, k, r+1}^{\{p\}} }
- \frac{ \sum_{j=i+1 }^{n-p }\gamma_{j, k, r}^{\{p\}} }{  \sum_{j=0 }^{n-p }\gamma_{j, k, r}^{\{p\}} }, \quad 0\leq i,r\leq n-p-1, 
\end{equation}
similar to (\ref{local_iterate}). The passage from level $p$ to level $p+1$ is successful when all quantities $\gamma_{i, k, r}^{\{p+1\}}$ are positive (at least in accordance with the positivity property of local expansions). The W-space $\E$ is an EC-space on $[a,b]$ if and only if the test is successful up to level $n-1$. In other words, if, for some integers $p, k, i$, $p\leq n-1$, some $0\leq k\leq q$, $0\leq i\leq n-p$, we obtain a negative $\gamma_{i, k, r}^{\{p\}}$, the  step stops, and the W-space $\E$ is not an EC-space.  

\smallskip
We have presented above the theoretical test on which the numerical test is based. In its numerical version, ``non-zero" or ``positive", etc,  is checked up to fixed tolerances, see \cite{BCM17}.

\begin{remark}
\label{test->poids}
If the test is successful up to level $(n-1)$, we know that it can be continued successfully until dimension one. In other words, not only can we say that $\E$ is an EC-space on $[a,b]$, but we can even add that
\begin{enumerate}[---]
\item with $w_{n-1}:=V_0^{\{n-1\}} +  V_1^{\{n-1\}}$ and $w_n:=V_0^{\{n\}}$, we have built  one system of weight functions (among infinitely many) associated with $\E$ --- \ie $\E=EC(w_0, \ldots, w_n)$ --- which  depends only on our choice of the initial Bernstein-like basis $(V_0^{\{0\}}, \ldots, V_n^{\{0\}})$; 
\item at each level, the  global Bernstein-like basis $(V_0^{\{p\}}, \ldots, V_{n-p}^{\{p\}})$ implicitly obtained in the test is a positive Bernstein-like basis relative to $(a, b)$ in the EC-space on $[a,b]$, $DL_{p-1}\E=EC(w_p, \ldots, w_n)$; 
\item at each level, dividing the equality $w_p= V_0^{\{p\}}+ \ldots+ V_{n-p}^{\{p\}}$ by $w_p$ provides us with  the Bernstein basis $(B_0^{\{p\}}, \ldots, B_{n-p}^{\{p\}})$ relative to $(a, b)$ in the space $L_p\E=EC(\un, w_{p+1}, \ldots, w_n)$; 
\item read in the reverse sense, the  relations between the bases at level $p$ and at level $p+1$ on which the test is based, produce  one among infinitely many recurrence formul\ae\  of the form (\ref{IRR}) \cite{JJA09}. 
\end{enumerate}
\end{remark}

\begin{example}
\label{HT_HTH}
Here we take $q=1$. Given $n\geq 1$, let $\E_0$ be the space spanned on $[t_0, t_1]$ by the functions $1, x, \ldots, x^{n-2}, \cos x, \sin x$ (to which we refer here as the {\it trigonometric space}), and let  $\E_1$ be the  space spanned on $[t_1, t_2]$ by the functions $1, x, \ldots, x^{n-2}, \cosh x, \sinh x$ (to which we refer here as the {\it hyperbolic space}). Setting $T:=t_1-t_0$ and $H:=t_2-t_1$, we want to know how to choose the pair $(T,H)$ so that the W-space $\E$ on $[t_0, t_2]$ obtained by the $C^n$ connection of $\E_0$ and $\E_1$ is an EC-space on $[t_0, t_2]$. To apply the test, we first have to make sure that  the space $\E_0$ is an EC-space on $[t_0, t_1]$. It is known that this is obtained by requiring $T<\ell_n$, where (see \cite{CMP04, CMP14, CMP17} and Subsection 4.3.1 below)
\begin{equation}
\label{ellT}
\ell_1=\pi, \quad \ell_2=\ell_3=2\pi, \quad \ell_4=\ell_5\approx 8. 9868, \quad \ell_6=\ell_7\approx11.5269, \quad \ell_8=\ell_9\approx13.9758, \quad\ldots
\end{equation}
On the left picture of Figure 1,  the region of the plane $(T,H)$ producing an EC-space on $[t_0, t_2]$ is located below the boundary curve, depending on the dimension. For  $n=1$ and $0<T<\pi$, this region,  limited by the green curve, coincides with the theoretical result obtained  in \cite{NUMA11}, that is,
$$
\E \hbox{ is an EC-space on }[t_0, t_2] \quad \Leftrightarrow \quad \cot T+\coth H >0. 
$$
We can see that the boundary curve presents more and more cusps as the dimension increases. The presence of such cusps can be explained by the various determinants (\ref{BLB?}) appearing in the first part of the test, each of them being expressed as a function of the two variables $T,H$. The different segments  of the boundary curve correspond to the different determinants (\ref{BLB?}) which effectively vanish at Step 0 in the test. 
The test could be applied with more sections as well. As an instance, in the right picture, we take three sections,  the first and last one are hyperbolic with length $H$, while the central one is  trigonometric with length $T$. For $n=8$, we design in the nine-dimensional space $\E$,  within the region of the plane $(T,H)$ ensuring that $D\E$ is an EC-space on $[t_0, t_3]=[0, 2H+T]$. with now $t_1$ and $t_2$ as the two interior knots. The curves are obtained with $H=2$ and, from top to bottom,  $T=0.1; 3.14;$ and  finally $T=6.6545$,  which corresponds to the limit curve.

\end{example}
  
\begin{figure}
\label{F1}
\begin{center}
\hskip-0.4cm
\includegraphics[width=7.8cm]{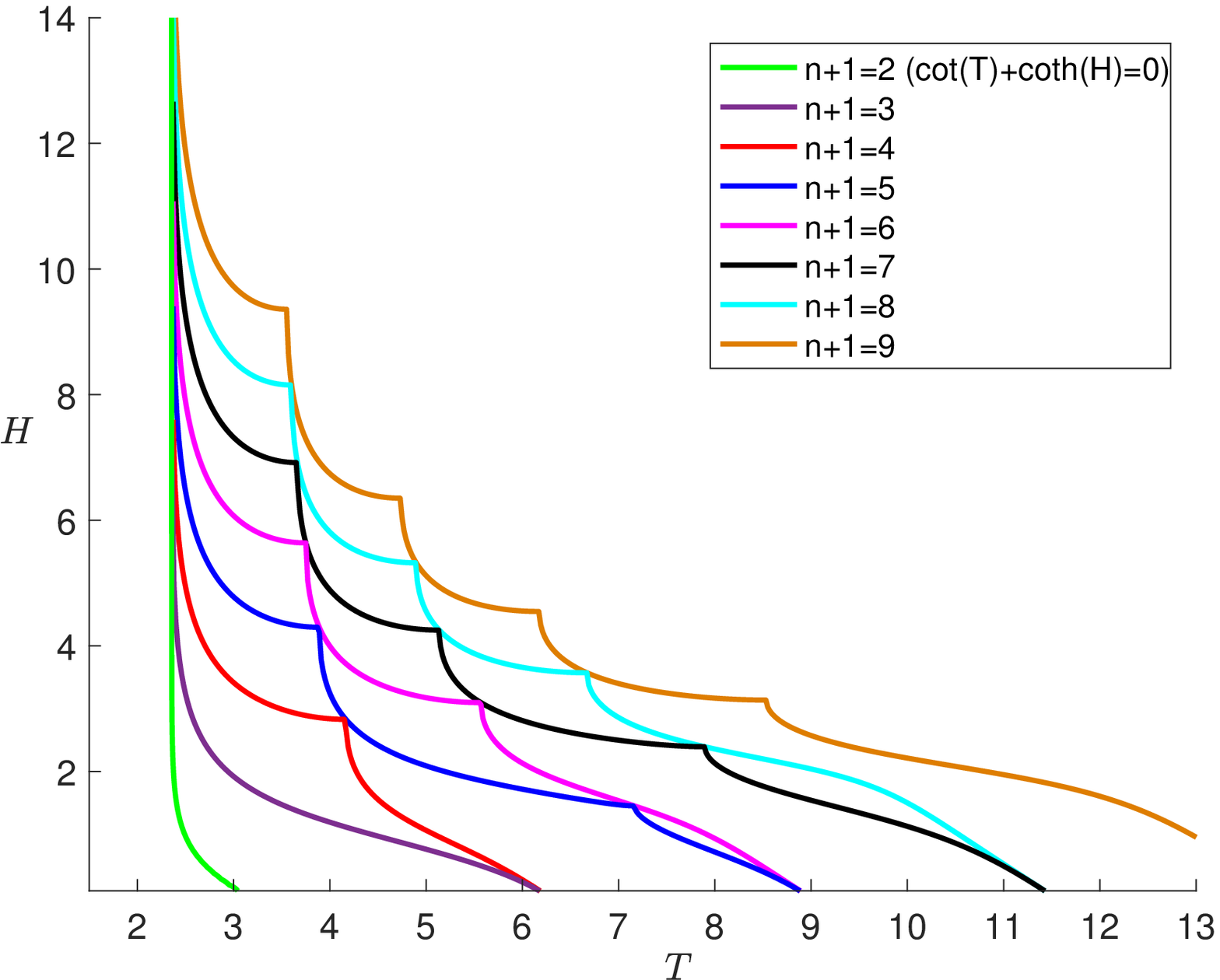}\hspace{1cm}
\raisebox{0.7cm}{\includegraphics[width=5cm]{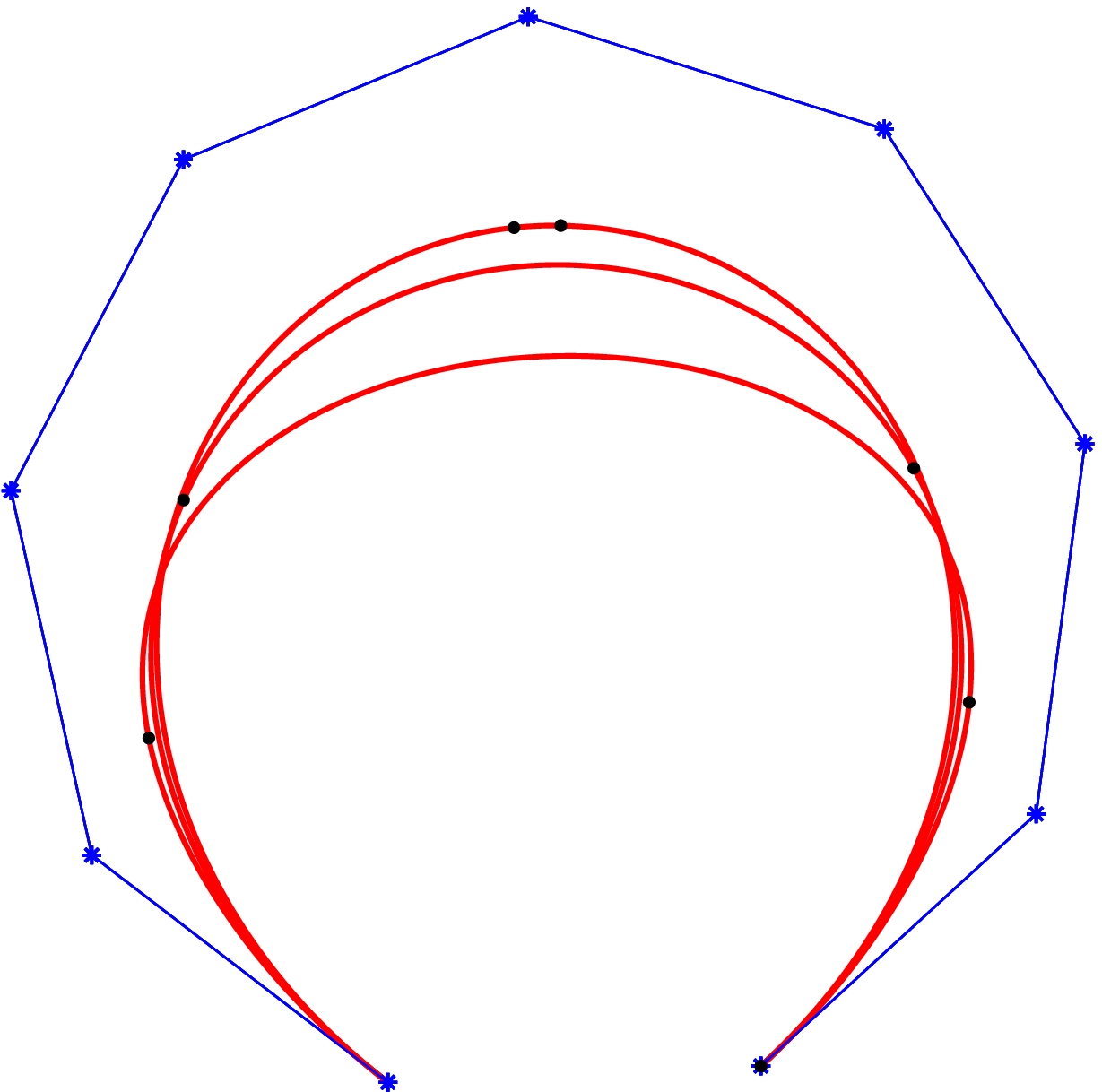}}
\end{center}
\caption{
{\bf Left:} In the plane $(T,H)$, for various values of $n$, the region producing a global EC-space on $[0, T+H]$ from a trigonometric space on $[0, T]$ and an hyperbolic space on $[T, T+H]$. {\bf Right:} A $C^8$ HTH-curve on $[0, 2H+T]$ for $H=2$ and $T=0.1; 3.14; 6.6545$ (limit curve), see Example \ref{HT_HTH}. }
\end{figure}

\section{Numerical Procedure for Critical Length: }

In this section we show how to use the numerical test briefly described in the previous one for computing the critical length of a given linear differential operator $L$ of order $(n+1)$ with constant coefficients, say 
\begin{equation}
\label{Ln}
L:=D^{n+1}+\sum_{i=0}^n a_iD^i, \quad \hbox{with } a_0, \ldots, a_n\in\RR, 
\end{equation}
where $n$ is a positive integer. Before explaining our approach, it is useful to review  some classical facts concerning such operators. 

\subsection{Preliminaries}
We are concerned with $\E_L=\ker(L)$. For a given $n$, the class of all such kernels coincides with the class of all $(n+1)$-dimensional W-spaces on $\RR$ which are closed under translation, or as well, with the class of all $(n+1)$-dimensional W-spaces on $\RR$ which are closed under differentiation. As is well known, the description of $\E_L$ follows from its characteristic polynomial 
\begin{equation}
\label{pn}
p_L(x):=x^{n+1}+\sum_{i=0}^n a_ix^i. 
\end{equation}
 In particular, $\E_L$  is closed under reflection (\ie for each $F\in \E_L$, the function $x\mapsto F(-x)$ belongs to $\E_L$) if and only if the polynomial $p_L$ is either even or odd. 

\begin{definition}
The critical length of the $L$ (or, as well, of $\E_L$) 
is  defined by
$$
\ell_L:=\sup\{h>0\ |\ \E_L \hbox{ is an EC-space on }[0, h]\}.
$$
\end{definition}

Due to $\E_L$ being  closed under translation, it is also the supremum of all $h>0$ such that $\E_L$ is an EC-space on any $[\alpha,\alpha+h]$. Two basic  facts are to keep in mind (see, for instance,  Proposition 1.15 in \cite{khaled}).   

\begin{proposition}
\label{firstprop}
The critical length $\ell_L$ lies in $]0, +\infty]$. More precisely
\begin{itemize}
\item if $p_L$ has only real roots, then $\ell_L=+\infty$;
\item otherwise,  $\frac{\pi}{M_L}\leq \ell_L<+\infty$, where $M_L$ stands for the maximum imaginary part of all roots of $p_L$, and $\E_L$ is not an EC-space on $[0, \ell_L]$.  
\end{itemize}
\end{proposition}

Finally, the classical way to compute  the critical length $\ell_L$ consists in determining the first positive zeros of a number of Wronskians, as  reminded below, \cite{CMP04, JJA12}. 

\begin{theorem}
Let $S$ be the unique element of $\E_L$ satisfying $S(0)=S'(0)=\dots=S^{(n-1)}(0)=0$,  $S^{(n)}(0)=1$. Then, the critical length $\ell_L$ can be obtained as
\begin{equation}
\label{comput_ellL}
	\begin{split}
	\ell_L&=\min_{0\leq k\leq n-1} \inf\{h>0\ |\ W(S,S', \ldots, S^{(k)})(h)\not=0\},\\
	         &=\min_{0\leq k\leq \frac{n-1}{2}} \inf\{\vert h\vert >0\ |\ W(S,S', \ldots, S^{(k)})(h)\not=0\}.
	\end{split}
\end{equation}
\end{theorem}
When all roots of $p_L$ are real, since $\ell_L=+\infty$, we can thus say that each of the Wronskians $W(S,S', \ldots, S^{(k)})$, $k=0, \ldots, n$, keeps the same strict sign on $]0, +\infty[$.
When  $\E_L$ is invariant under reflection, formul\ae\  (\ref{comput_ellL})  can be simplified as follows:
\begin{corollary}
Suppose that  $p_L$ is either odd or even.  Then we have
\begin{equation}
\label{comput_ellL(pair)}
	\ell_L=\min_{0\leq k\leq \frac{n-1}{2}} \inf\{ h>0\ |\ W(S,S', \ldots, S^{(k)})(h)\not=0\}.
\end{equation}
\end{corollary}

Consider the operator $\widehat L:=D\circ L$. Given that $\E_L=D\E_{\widehat L}$, according to Definition \ref{def:GFD}, we can say that  
$$
\ell_L:=\sup\{h>0\ |\ \E_{\widehat L} \hbox{ is an EC-space good for design on }[0, h]\}. 
$$
This justifies the following terminology. 
\begin{definition}
The critical length $\ell_L$ is called {\it the critical length for design of the operator $\widehat L:=D\circ L$}. 
\end{definition}
\noindent
The space $\E_{\widehat L}$ being obtained from $\E_L$ by integration, it should be observed that its critical length for design  is always less than or equal to its critical length.

\subsection{An alternative approach}

In this subsection, we assume that $p_L$ has at least one non-real root, and as previously, we denote by $M_L$ the maximum imaginary part of all non-real roots of $p_L$. We will compute the critical length $\ell_L\in ]0, +\infty[$ via the test recalled in the previous section. The computation comprises two successive parts:

\subsubsection{Step 1: Rough estimate of $\ell_L$}

We select a positive number $\ell_0<\frac{\pi}{M_L}$, close to $\frac{\pi}{M_L}$. From Proposition \ref{firstprop}, we know that $\E_L$ is an EC-space on $[0, \ell_0]$. Accordingly, we can state that
$$
\hbox{there exists an integer }\mu\geq 1 \hbox{ such that } \mu\ell_0<\ell_L\leq (\mu+1)\ell_0. 
$$
This integer is equivalently determined by the two properties below
\begin{equation}
\label{mu}
\E_L\hbox{ is an EC-space on }[0, \mu\ell_0] \hbox{ and } \E_L\hbox{ is not an EC-space on }[0, (\mu+1)\ell_0].
\end{equation}
Observing that, for each non-negative $k$, the restriction $\E_k$ of $\E_L$ to $[t_k, t_{k+1}]:=[k\ell_0, (k+1)\ell_0]$ is an EC-space on $[t_k, t_{k+1}]$, we can apply the test on $[t_0, t_{k+1}]$, with the sequence $\T=(t_1, \ldots, t_k)$ of interior knots, to determine whether $\E_L$ is an EC-space  on $[0, (k+1)\ell_0]$, successively for $k= 1, 2, \ldots$. The integer $\mu$ satisfying (\ref{mu}) is the first integer $k$ for which we obtain a negative answer. 

\subsubsection{Step 2: Search for $\ell_L$ in $]\mu\ell_0,  (\mu+1)\ell_0]$}
The second step consists in localising $\ell_L$ within the interval $]\mu\ell_0,  (\mu+1)\ell_0]$ by dichotomy. On account of (\ref{mu}), 
we first want to test if $\E_L$ is an EC-space on $[0, \mu\ell_0 +\frac{\ell_0}{2}]$. Then, 
\begin{enumerate}[--]
\item if the answer is affirmative, test if $\E_L$ is an EC-space on $[0, \mu\ell_0 +\frac{\ell_0}{2}+ \frac{\ell_0}{4}]$; 
\item if the answer is negative, test if $\E_L$ is an EC-space on $[0, \mu\ell_0 +\frac{\ell_0}{2}- \frac{\ell_0}{4}]$; 
\end{enumerate}
Continue the same way, that is, at each step increment or decrement the interval length by  $\frac{\ell_0}{2^n}$, until $\frac{\ell_0}{2^n}$ is less than a given tolerance. 

At each dichotomy step, we apply the test of Section 3. However,  in order to avoid numerical problems, we do not apply it using the same first $\mu$ intervals of length $\ell_0$ and adding an interval which might be of smaller and smaller length,  but using only two consecutive intervals both of same variable length. To be more precise, in the first step we  take 
$$
[t_0, t_2]:=\bigl[0, \ell_0 \bigl(\mu+\frac{1}{2}\bigr)\bigr], \quad t_1=\frac{t_2}{2}=\frac{\ell_0}{2} \bigl(\mu+\frac{1}{2}\bigr), \quad \E_0:={\E_L}_{\big\vert[t_0, t_1]}, \quad \E_1:={\E_L}_{\big\vert[t_1, t_2]}.
$$
Since $\mu\geq 1$, for $i=0, 1$, $\E_i$ is an EC-space on $[t_i, t_{i+1}]$, which enables us to the apply the test. If the answer is affirmative (resp., negative) we do the same, after incrementing (resp., decrementing) $t_2$ by $\frac{\ell_0}{4}$, and so forth. At each step of the dichotomic process, each  $\E_i$, $i=0, 1$,   is guaranteed to be an EC-space on $[t_i, t_{i+1}]$.

\subsection{Examples}
We illustrate the previous procedure with several  instances of differential operators, indicated by their characteristic polynomials. 
We have chosen to illustrate the behaviour of the critical lengths up to dimension nine,
since this has proven to be sufficient to formulate theoretical conjectures and, at the
same time, it comprises all cases of practical interest for applications (\eg design).
The computational method of the critical length has no limitations on the dimension of
the space, but being a numerical method and working at the limits of the critical lengths
it is subject to the machine precision and to the fixed tolerances.
The following results are obtained in Matlab with tolerances fixed to 1.0e-30 for the test of Section 3, and 1.0e-10 for the dichotomy procedure. 

\subsubsection{Basic example: $p_n(x)=x^{n-1}(x^2+b^2)$, $b>0$}
 This is the simplest class of kernels $\E_n$ of differential operators for which the critical lengths are not infinite,  spanned  on $\RR$ by the $(n+1)$ functions $1, x, \ldots, x^{n-2}, \cos(bx), \sin(bx)$ (cycloidal or trigonometric spaces depending on the literature). The corresponding critical lengths $\ell_n(b)$ are shown in Figure 2, top left (under the name $\ell_L$), as functions of the variable $b$, for increasing values of  $n$. 
Clearly, we have 
$$
\ell_n(b)=\frac{\ell_n(1)}{b}, \quad b>0.
$$ 
This  explains why each curve resulting from our numerical procedure is a branch of hyperbola. 
Moreover, the first values of $\ell_n:=\ell_n(1)$ are those already given in (\ref{ellT}). The critical length $\ell_n$ was studied  in several papers by Carnicer, Mainar, Pe\~na \cite{CMP04, CMP14, CMP17}, their main results being stated below: 

\begin{theorem}
\label{th:T}
The successive critical lengths $\ell_n$ of the spaces $\E_n$ spanned by $1, x, \ldots, x^{n-2}, \cos x, \sin x$, $n\geq 1$,  satisfy
\begin{equation}
\label{eq:T}
\ell_{2k}=\ell_{2k+1}=2 j_{k-\frac{1}{2}, 1}<\ell_{2k+2}, \quad k\geq 1,  
\end{equation}
where, for each positive $\alpha$, $j_{\alpha, 1}$ stands for the first positive zero of the Bessel function of the first kind:
$$
J_\alpha(x)=\sum_{m=0}^\infty\frac{(-1)^m}{m! \ \Gamma(m+\alpha+1)}\left(\frac{x}{2}\right)^{2m+\alpha}\!\!.
$$ 
\end{theorem}

\noindent
On the other hand, since $p_n$ is either odd or even, we know that the critical lengths can be computed by (\ref{comput_ellL(pair)}). It is worthwhile mentioning that, whatever the dimension, 
$$
\ell_n\hbox{  is the first positive zero of
the Wronskian }W(S_n, \ldots, S_n^{(p)}),  \hbox{  where }p:=\left\lfloor \frac{n-1}{2}\right\rfloor, 
$$
where the non-zero function $S_n\in \E_n$ vanishes $n$ times at 0, and where $\lfloor \pt \rfloor$ stands for the floor function. As shown in  \cite{CMP14}, none of the other Wronkians involved in (\ref{comput_ellL(pair)}) vanishes on $]0, +\infty[$.

Our procedure, first ran in this class of spaces, is in perfect accordance with Theorem \ref{th:T}. In particular, the two-by-two behaviour can easily be observed in Figure 2.  

\subsubsection{Other examples of the form $p_n(x)=x^{n-3}p_3(x)=x^{n-3}(x^4+a_2x^2+a_0)$, $n\geq 3$}
Let us recall that the  four-dimensional spaces $\E_3$ associated with characteristic polynomials $p_3(x)=x^4+a_2x^2+a_0$, with at least two non-real roots,  were thoroughly investigated in \cite{BM12}.  It is natural to apply the procedure  in the class of spaces   obtained from them by repeated integration, by comparison with its simplest subclass of cycloidal spaces.  In spite of this class being relatively limited and simple, we will already observe different  interesting behaviours depending on the polyniomal $p_3$ and on the integer $n$.  This will clearly point out the difficulty of determining  the critical length for any given differential operator, and the impossibility to easily foresee what it will be. 

\medskip
\noindent
$\bullet$ $p_3(x)=(x^2-a^2)(x^2+b^2)$, $a, b>0$:\par
\smallskip

The space $\E_n$ is spanned by the $(n+1)$ functions $1, x, \ldots, x^{n-4}, \cosh(ax), \sinh(ax), \cos(bx), \sin(bx)$. The critical length  of $\E_n$ depends on the two parameters $a, b$, and we denote it by $\ell_n(a,b)$. Now, we clearly have
\begin{equation}
\label{mult_alpha}
\ell_n(\alpha a,\alpha b)=\frac{1}{\alpha}\, \ell_n(a,b), \quad \alpha, a,b>0. 
\end{equation}
This makes it sufficient to apply our numerical procedure for the computation of  $\ell_L:=\ell_n(1,b)$. The corresponding graph   in function of the single variable  $b$  is shown for increasing values of $n\geq 3$ in Figure 2, top right. For $n=3$, the procedure confirms the results obtained  in \cite{BM12}, namely the fact that the critical length $\ell_3(a,b)$  is the only solution of the equation
\begin{equation}
\label{ZH3}
(b^2-a^2)\sinh(ax)\sin(bx)=2ab\bigl(1-\cosh(ax)\cos(bx) \bigr), \quad x\in \left]{\pi\over b}, {2\pi\over b}\right[. 
\end{equation}
With the same notation as for the the cycloidal example, it  is the first positive zero of the Wronskian $W(S_3, S_3')$, while the function $S_3$  has no positive zero.  From Figure 2, we conjecture  that this class of spaces follows the same two-by-two behaviour as the cycloidal spaces (apart from the connection with Bessel functions) involving only one of the Wronskians in formula (\ref{comput_ellL(pair)}). Nevertheless, our purpose here is not  to solve this conjecture.  \par
\begin{figure}
\includegraphics[width=7.7cm]{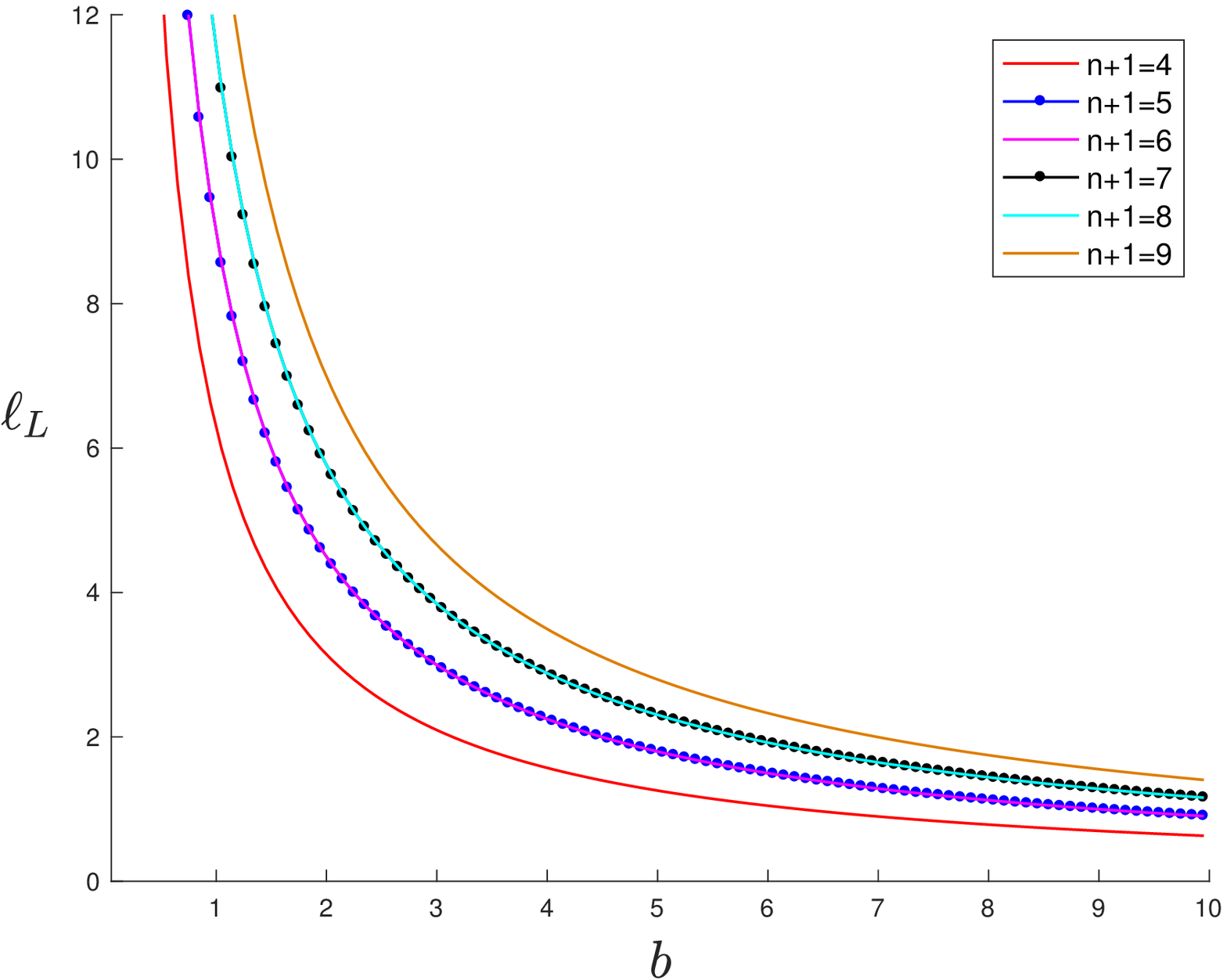}
\includegraphics[width=7.7cm]{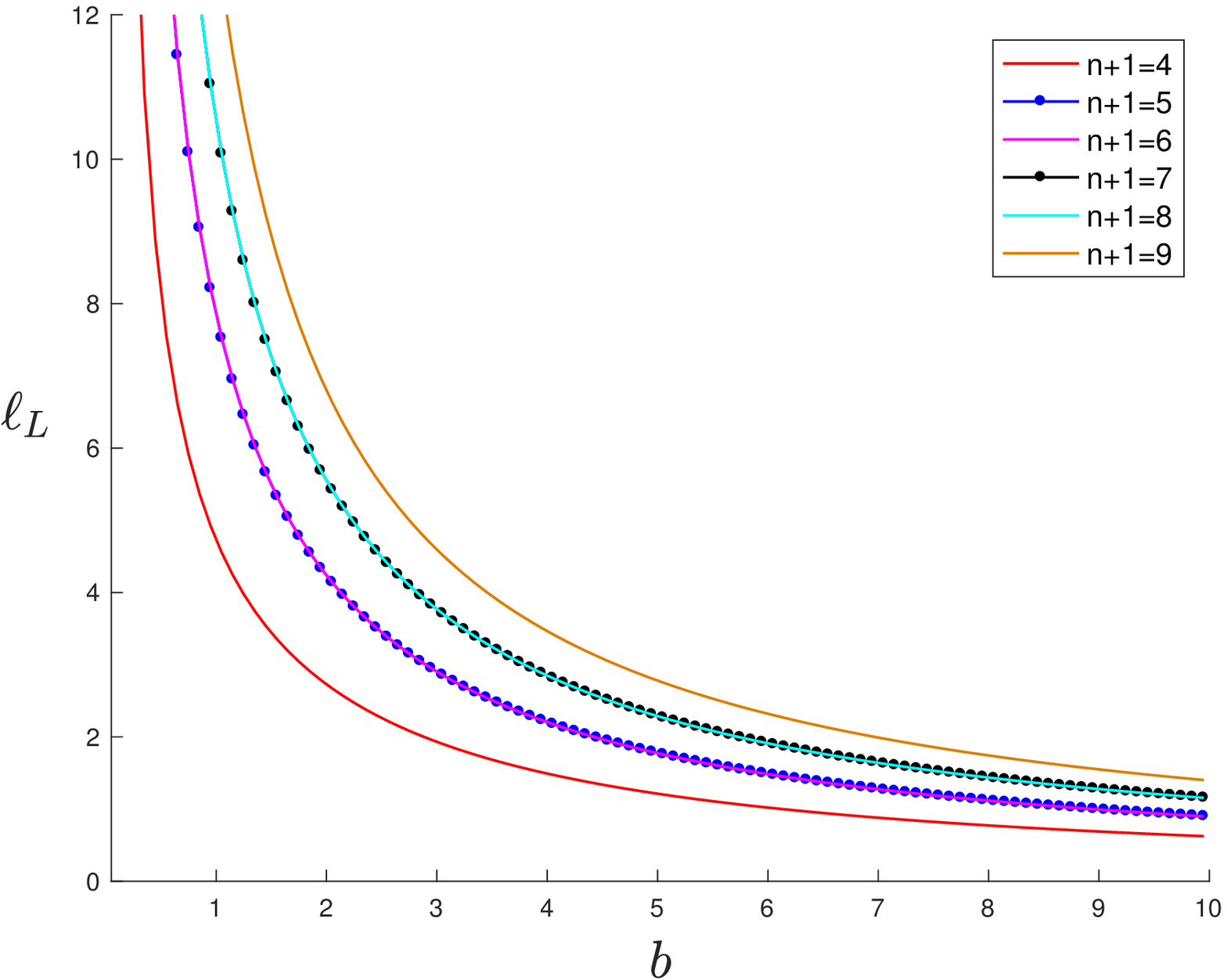}
\includegraphics[width=7.7cm]{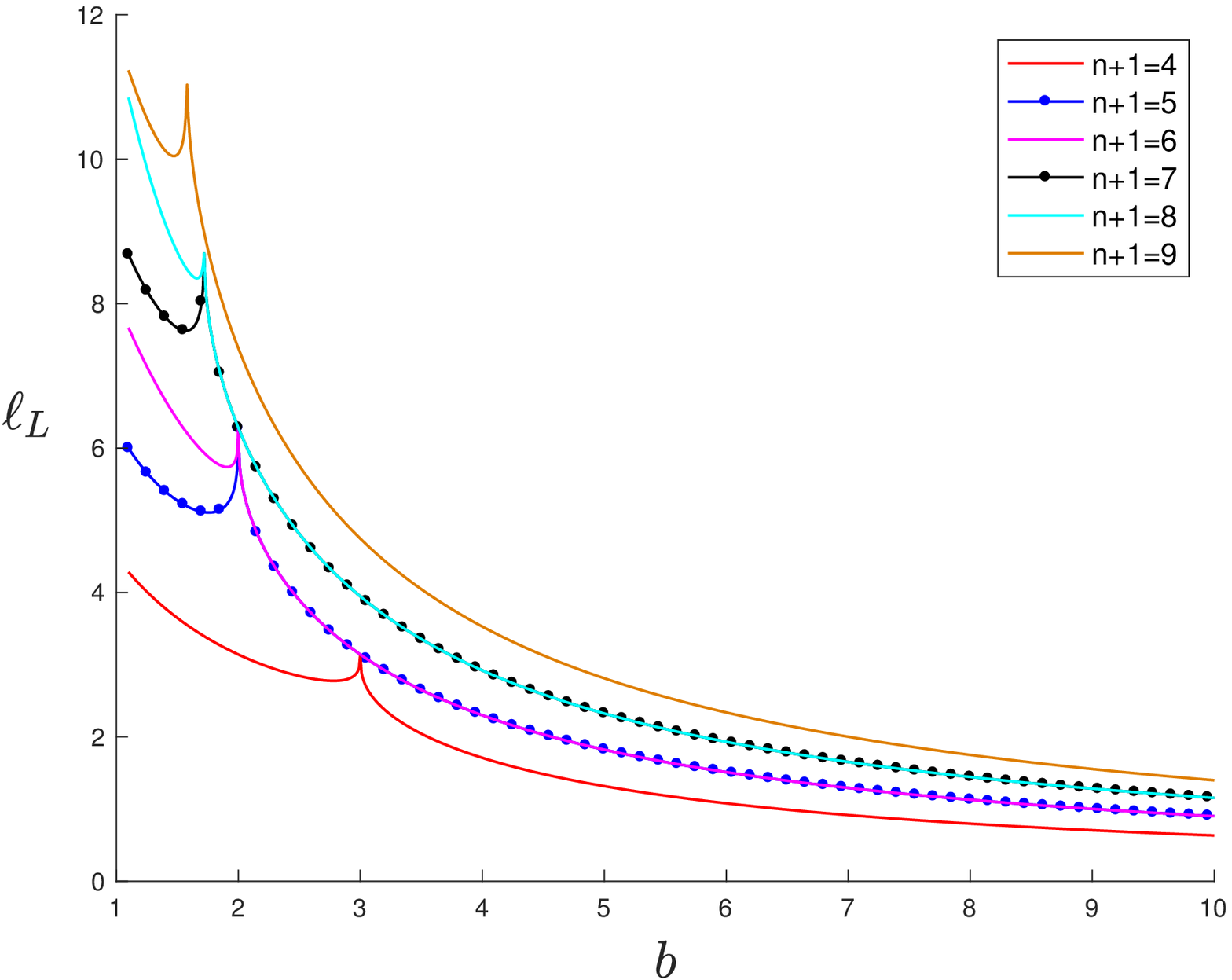}
\includegraphics[width=7.7cm]{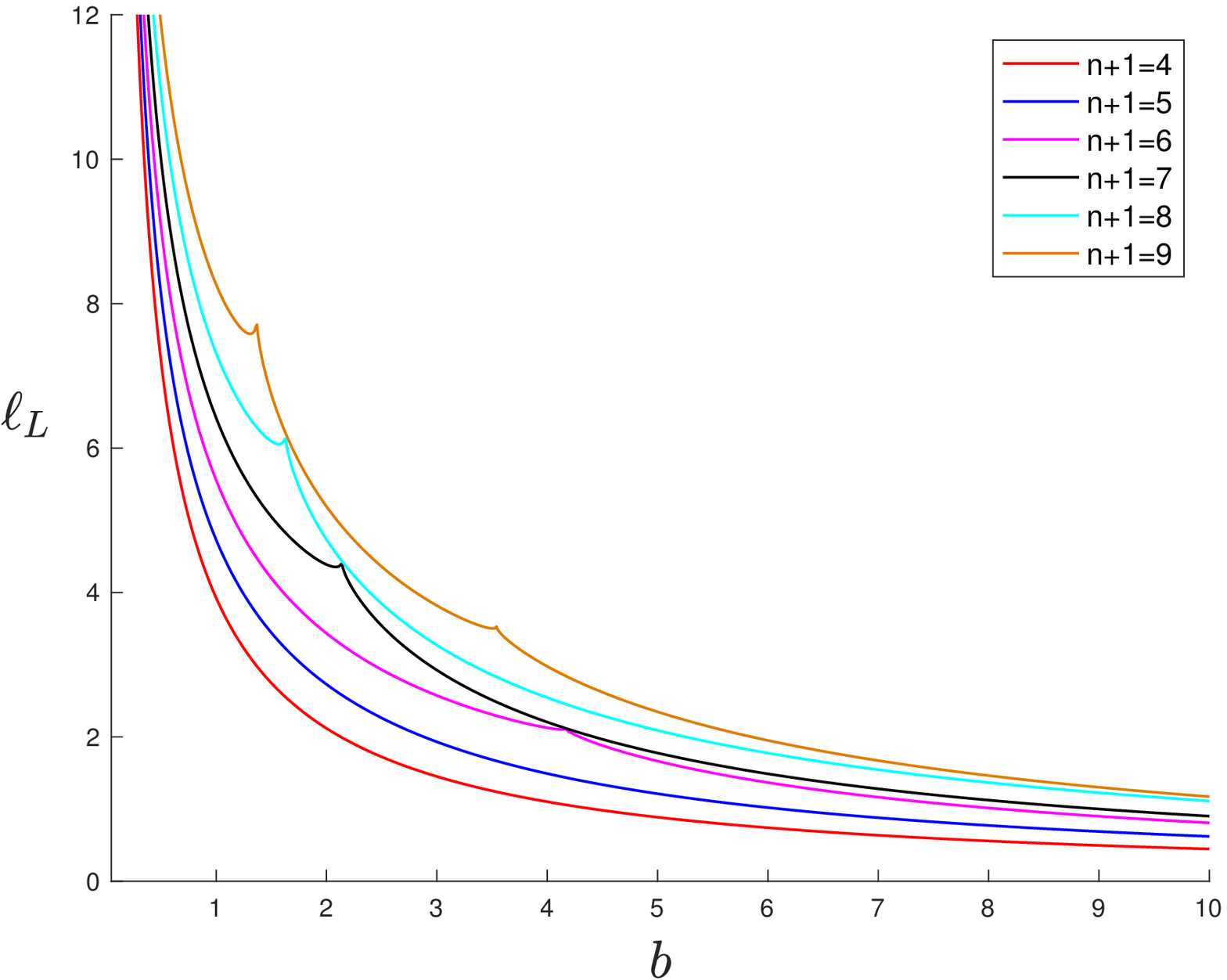}
\caption{
The critical length of the space $\E_L$ as a function of the positive parameter $b$, where $\E_L$ is spanned by $1,x, \ldots, x^{n-2}, \cos (bx), \sin (bx)$  (top left); $1,x, \ldots, x^{n-4}, \cosh x, \sinh x, \cos (bx), \sin (bx)$ (top right); $1,x, \ldots, x^{n-4}, \cos x, \sin x, \cos (bx), \sin (bx)$, with $b>1$ (bottom left); $1,x, \ldots, x^{n-4}$, $\cosh x\cos(bx), \cosh x\sin(bx), \sinh x\cos (bx), \sinh x\sin (bx)$ (bottom  right).}
\end{figure}
\eject\noindent
$\bullet$ $p_3(x)=(x^2+a^2)(x^2+b^2)$, $a,b>0$:\par
\smallskip

We first assume that $a\not=b$. Here, the space $\E_n$ is spanned by the $(n+1)$ functions $1, x, \ldots, x^{n-4}$, $\cos(ax), \sin(ax), \cos(bx)$, $\sin(bx)$. Assuming that $a<b$, and taking account of (\ref{mult_alpha}), we apply the numerical procedure to compute $\ell_L:=\ell_n(1,b)$, for $b>1$. For each $n\geq 3$, the presence of one cusp in the graph indicates that the two parts of the graph are  related to two different Wronskians. Second interesting observation: the right parts of the graphs follow the same two-by-two behaviour as the cycloidal spaces, while the first parts change whenever $n$ increases. The results obtained in \cite{BM12} for $n=3$ are confirmed, namely: 
\begin{enumerate}[--]
\item for $a<b\leq 3a$, $\ell_3(a,b)$ is the first zero of $S_3$, and it is the only solution of the equation
$$
b(\sin(ax)=a\sin(bx), \quad x\in\left[\,\frac{\pi}{b}\left\lfloor\frac{b}{a}\right\rfloor\,, \,\frac{\pi}{b}\left\lceil\frac{b}{a}\right\rceil\,\right], 
$$
where  $\lceil\pt\rceil$ denotes the ceiling function; 
\item for $b\geq 3a$, $\ell_3(a,b)$ is the first zero of $W(S_3, S_3')$, and it is the only solution of the equation
$$
(b-a)\sin\left(\frac{(b+a)x}{2}\right)+(b+a)\sin\left(\frac{(b-a)x}{2}\right)=0, \quad x\in \left]\frac{2\pi}{b}, \frac{2\pi}{b-a}\right[.
$$
\end{enumerate}

When $a=b$, the space $\E_n$ is spanned by the $(n+1)$ functions $1, x, \ldots, x^{n-4}, \cos(bx), \sin(bx)$, $x\cos(bx)$, $x\sin(bx)$, and the graphs of the critical lengths $\ell_n(b)=\ell_n(1)/b$  (not shown here) change whenever the dimension increases, with $\ell_n(1)$ in accordance with the left parts of the case $a=1<b$ of which it is the limit situation. \par
\medskip
\noindent
$\bullet$ $p_3(x)=x^4+2(b^2-a^2)x^2+(a^2+b^2)^2$, $a,b>0$:\par
\smallskip

Here, the space $\E_n$ is spanned on $\RR$ by the $(n+1)$ functions $1, x, \ldots, x^{n-4}, \cosh(ax)\cos(bx)$, $\sinh(ax)\cos(bx)$, $\cosh(ax)\sin(bx), \sinh(ax)\sin(bx)$. 
In \cite{BM12}, the critical length $\ell_3(a,b)$ was proved to be  the only solution of the equation 
\begin{equation}
\label{ZS9}
b\tanh(ax)=a\tan(bx), \quad x\in \left]{\pi\over b}, {3\pi\over 2 b}\right[.
\end{equation}
It  is  the first positive zero of $S_3$, while the Wronskian $W(S_3, S_3')$ has no positive zero. As $n$ increases,  we can see, in the graph of $\ell_L:=\ell_n(1,b)$,  the presence of one cusp (for $n=5,6,7$), then two cusps ($n=8$),  indicating that two (resp. three) of the Wronskians appearing in formula (\ref{comput_ellL(pair)}) are involved in the graph.

\subsubsection{More examples}

These additional examples are intended to show that our numerical procedure is not limited to classes of spaces obtained by repeated integration. \par

\medskip
\noindent
{$\bullet$ More pairs of trigonometric functions:}

\smallskip
The critical length of the $(2n)$-dimensional space $\E_{2n-1}$ spanned by the functions $\cos(x), \sin (x)$, $\cos(2x), \sin (2x), \ldots,  \cos(nx), \sin (nx)$ is equal to $\pi$. This well-known result (see, \eg \cite{schu}) was successfully checked through the numerical procedure for several values of $n$. But what about other pairs of trigonometric functions? A very partial answer to this question is given for   three pairs $\cos(ax), \sin (ax)$, $\cos(bx), \sin (bx), \cos(cx)$, $\sin (cx)$, with $0<a<b<c$. Without loss of generality, we take $a=1$, and for several values of $b>1$, we show in Figure 3 the graph of the critical length as a function of the only variable $c>b$. As a generalisation of what happened with the four-dimensional space spanned by two different pairs (see Figure 2, bottom left) we can see in these examples that the curve presents two cusps (\ie three active Wronskians) and we can observe how they evolve as $b$ increases. However, these few pictures are certainly not enough to state any possible conjecture.

\medskip
\noindent
{$\bullet$ A non-symmetric example:}

\smallskip

Of course, the procedure is not limited to spaces closed under reflection.  This is the reason why we conclude the illustrations with the characteristic polynomial
$p_4(x) = x^2(x-a)(x^2+b^2)$, where $a$ is any non-zero real, and where $b>0$.  In other words, $\E_4$ is spanned by $1, x, e^{ax}, \cos(bx), \sin(bx)$. In Figure 4,  we show  the graph of the critical length $\ell_4(a, 1)$ as a function of $a\not=0$, and the graph of $\ell_4(1, b)$ as a function of $b>0$. Note that the limit value $\ell_4(0, 1):=\lim_{a\to 0}\ell_4(a, 1)$ is the value $\ell_4=8.9868$ given in (\ref{ellT}).

 \begin{figure}
\begin{center}
\includegraphics[height=0.25\textwidth]{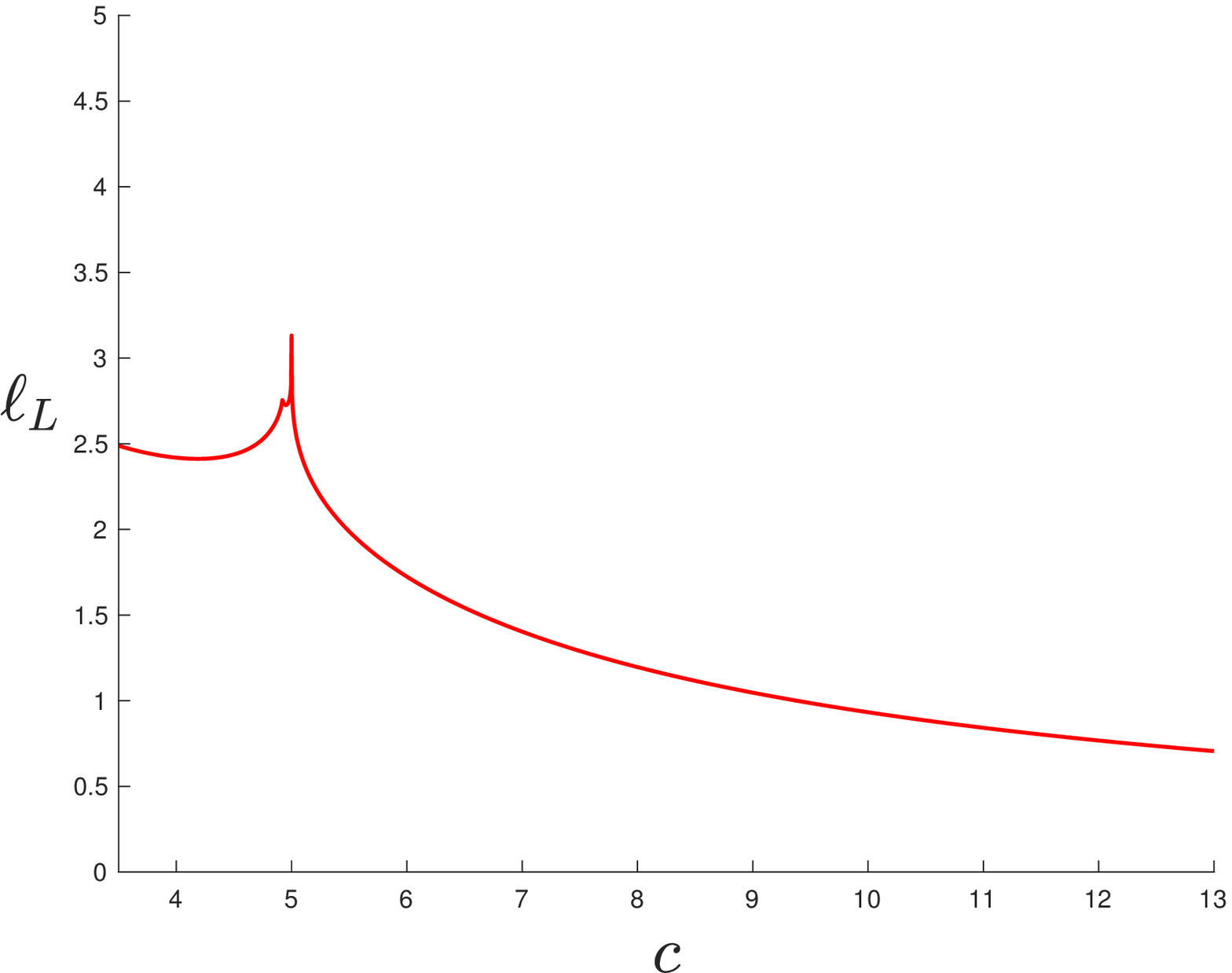} \hfill
\includegraphics[height=0.25\textwidth]{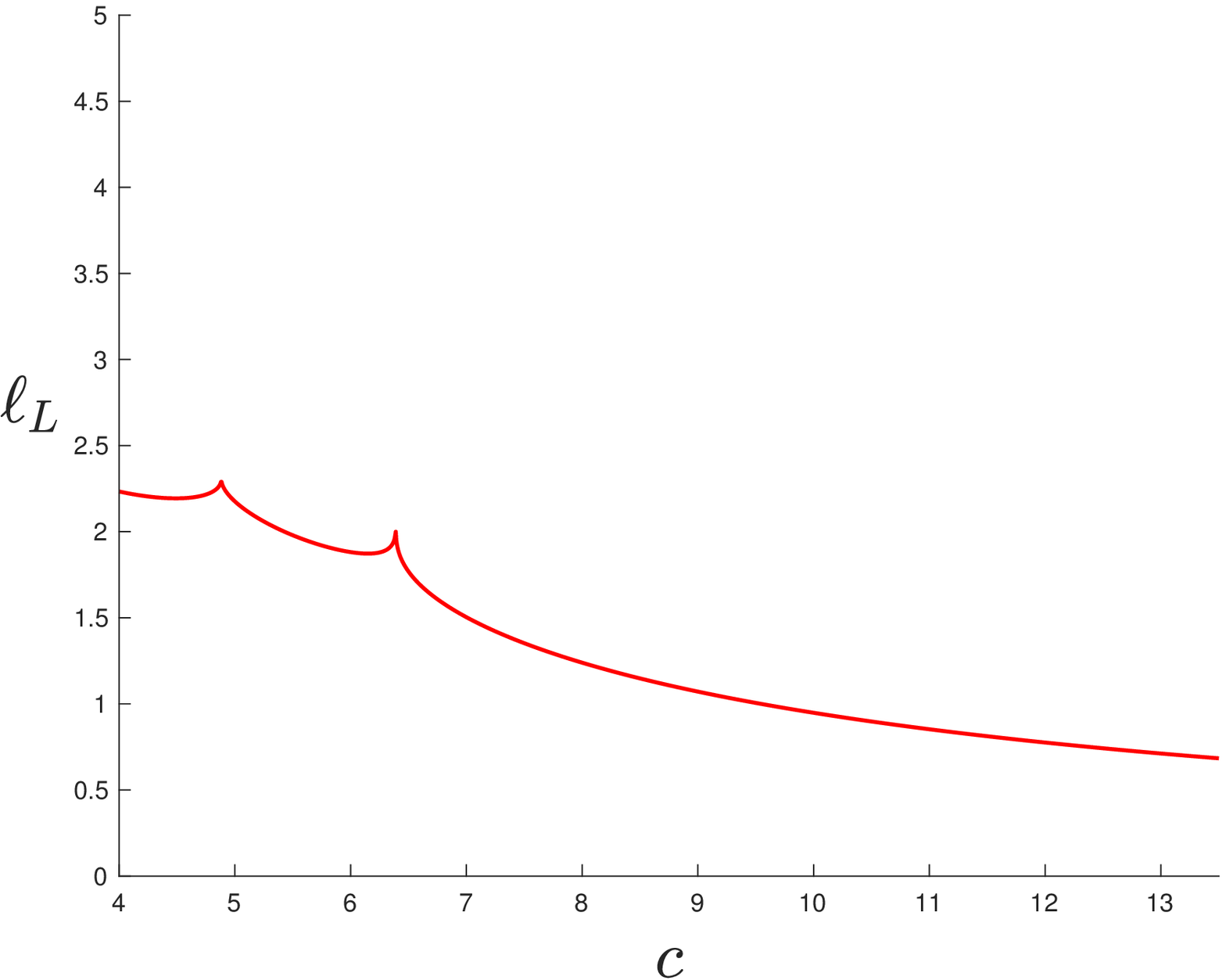} \hfill
\includegraphics[height=0.25\textwidth]{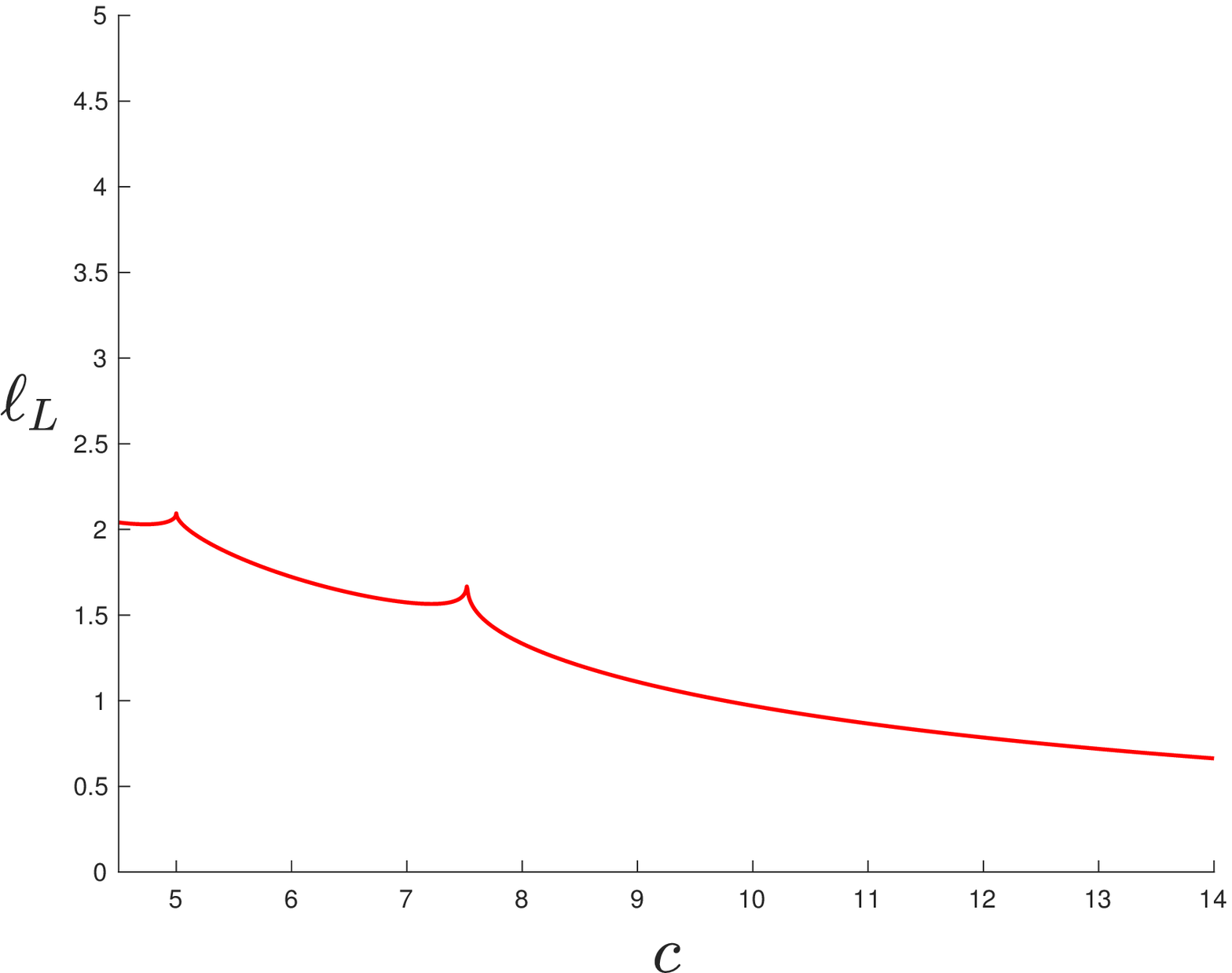} \hfill
\tikzstyle{every picture}+=[remember picture,overlay]
\begin{tikzpicture}
\node[circle,draw=black,minimum width=0.45cm] at (-14.15,2.32) {};
\node[anchor=center] at (-12.5,2.2){\includegraphics[width=0.45cm]{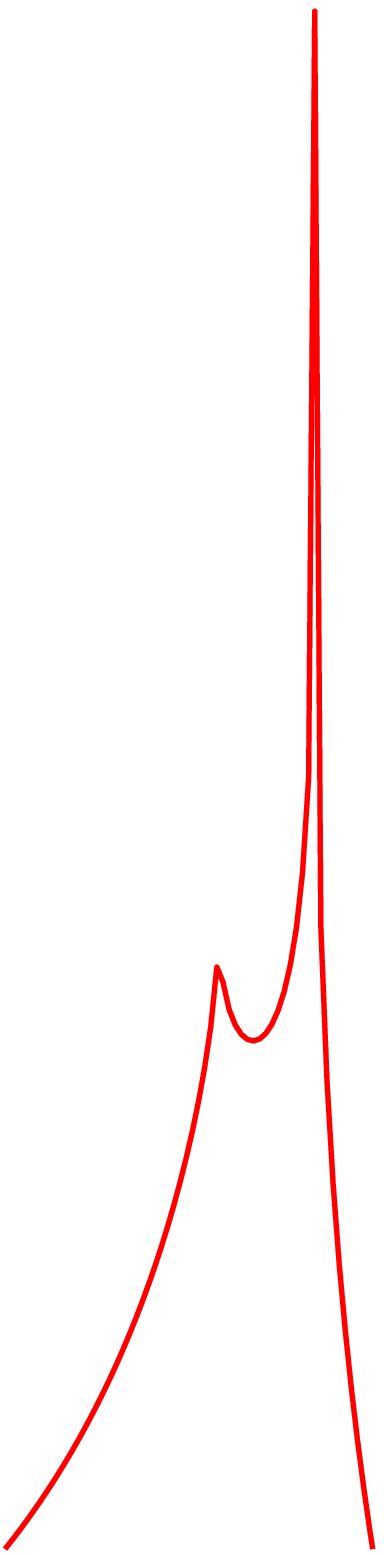}};
\node[circle,draw=black,minimum width=1.9cm] at (-12.5,2.2) {};
\end{tikzpicture}
\caption{Critical length for the space $\E_n$ spanned by $\cos ax, \sin ax,\cos bx, \sin bx,\cos cx, \sin cx,$ with $a=1$ and, from left to right, $b=3,3.5,4$.}
\end{center}
\begin{center}
\includegraphics[height=0.25\textwidth]{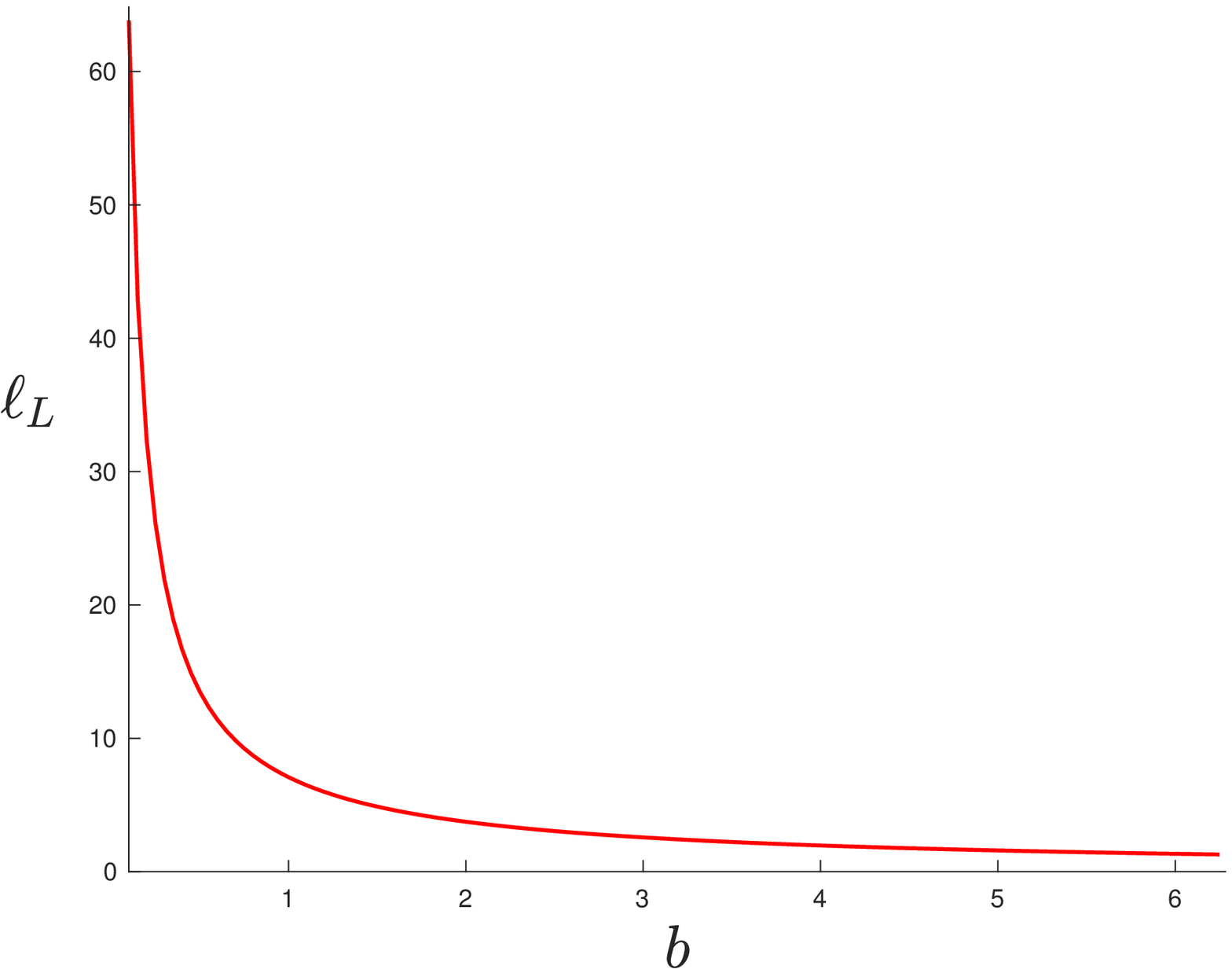} \hspace{1cm}
\includegraphics[height=0.25\textwidth]{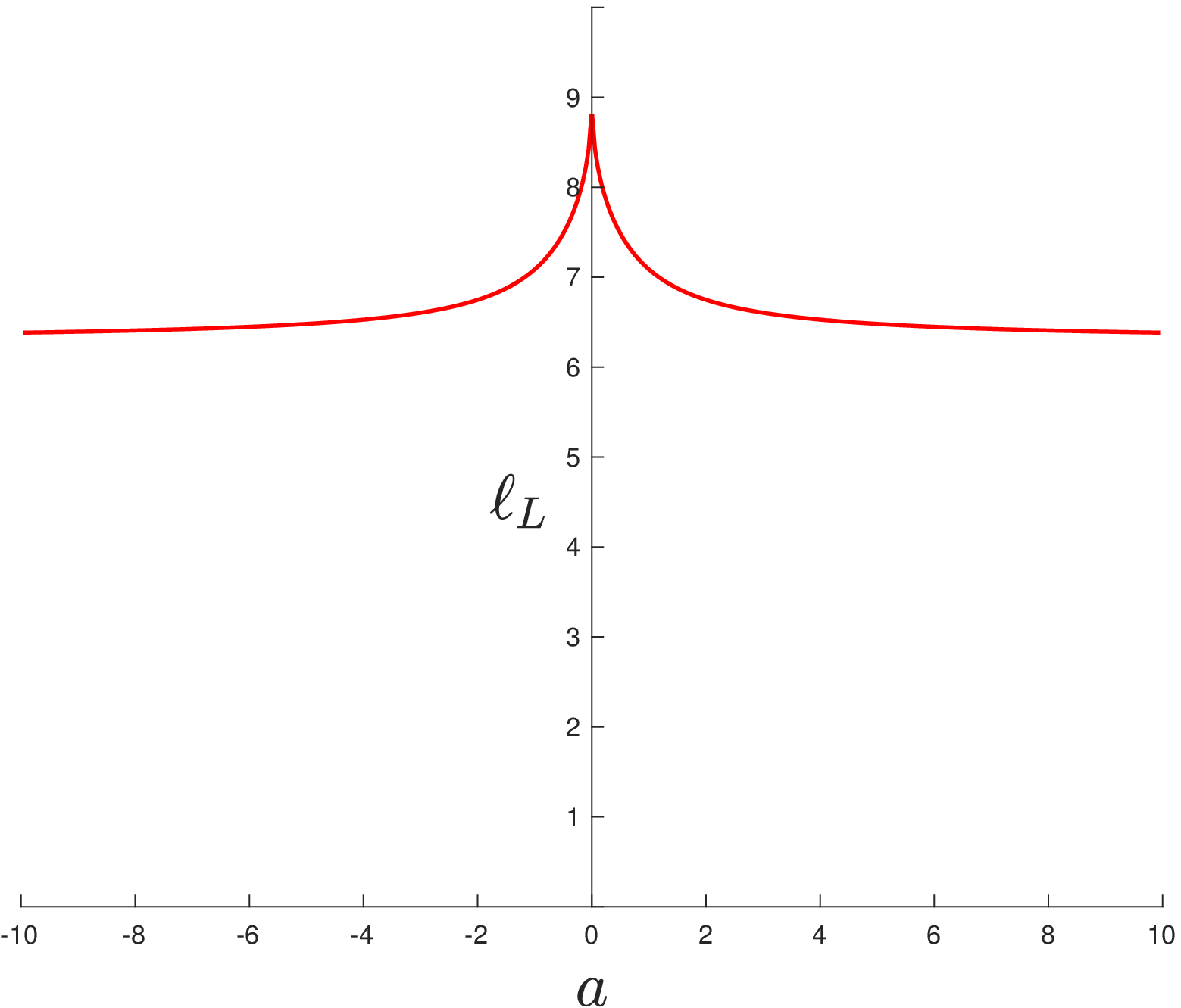} 
\caption{Critical length for the space $\E_4$ spanned by $1,x,e^{ax}, \cos bx, \sin bx$ with $a=1$ (left) and $b=1$ (right).}

\end{center}
\end{figure}

\section{Critical length: why?}

In this concluding section, we gather crucial remarks which clearly point out the importance of knowing the critical lengths, and thus the interest of our numerical method to compute them. 

\medskip

Due to Theorem \ref{th:EC(w)}, it may seem tempting to consider systems of weight functions as the starting point of anything concerning EC-spaces on a given closed bounded interval $[a,b]$,  all the more so as the classical approximation results generalising the polynomial framework are based on weight functions \cite{karlin, schu}.  If one adopts  this angle, the Bernstein bases relative to $[a,b]$ are then naturally provided by the associated integral recurrence relations (\ref{IRR}). This approach, thoroughly developed in \cite{bister, bisterP}\footnote{These references address the general framework of integral-positive weight functions, with application to EC-spaces and associated splines. Integral recurrence relations can also be understood through the blossoming approach, as shown in \cite{JJA09}, including EC-piecewise spaces or Quasi EC-spaces, and associated splines, see also  \cite{howto, NUMA11}. } was certainly worth attention for theoretical purposes when it was elaborated.  It was also commonly used to develop examples,  in particular with  ordinary (rather than generalised) integration, see for instance \cite{ChenWang, MP, wang}. 
It is important though to draw the reader's attention to its limitations on the practical side. 

First of all, if we start  with any arbitrary system $(w_0, \ldots, w_n)$ of weights functions on $[a,b]$, chances are that the  associated space $\E_n:=EC(w_0,w_1, \ldots, w_n)$ will not present much interest for design or any other application. As a general rule, replacing polynomial spaces by given EC-spaces is motivated by either geometric constraints or likely to happen shape effects, or whatever purpose, but always with the space $\E_n$ of functions which we expect to work in as the starting point.  Unless the dimension is really small, neither on which intervals it is an EC-space (possibly good for design) nor possible associated systems of weight functions are easy-to-solve questions.  This is why constructing Bernstein bases through (\ref{IRR}) may be extremely limiting.  With a view to stress this limitation, we take the simplest class of trigonometric spaces spanned by the function $1, x, \ldots, x^{n-2}, \cos x, \sin x$, with the following notations:
$$
L_n=D^{n+1}+D^{n-1}, \quad \ell_n:=\hbox{critical length of }L_n, \quad n\geq 1. 
$$
Subsequently, we work with a fixed integer $n\geq 1$, and a fixed interval $[a,b]$, with $0<b-a<\ell_n$. For each $k\geq 0$, the notations $\E_k$, $\PP_k$, stand for the restrictions to $[a,b]$ of $\ker L_k$, and of the degree $k$ polynomial space, respectively. From the condition $b-a<\ell_n$ and from Theorems \ref{th:EC(w)} and \ref{th:w0}, we know that one can find infinitely many nested sequences  
\begin{equation}
\label{Estar}
\E_0^\star\subset \E_1^\star\subset \cdots\subset \E_{n-1}^\star\subset \E_{n}^\star:=\E_n, 
\end{equation}
where, for $i=0, \ldots, n-1$, $\E_i^\star$ is an $(i+1)$-dimensional  EC-space (or W-space as well) on $[a,b]$, and the search for such nested sequences is equivalent to the search for systems $(w_0, \ldots, w_n)$ of weights functions on $[a,b]$ such that $\E_n=EC(w_0, \ldots, w_n)$. Now, the space $\E_n$ contains two obvious nested sequences of W-spaces
\begin{equation}
\label{twonested}
\PP_0\subset\PP_1\subset \cdots\subset \PP_{n-2}\subset \E_n, 
\qquad \E_1\subset \E_2\subset \cdots\subset \E_{n-1}\subset \E_n. 
\end{equation}
Each of them is ``almost" a sequence (\ref{Estar}), and we will subsequently discuss when it is ``exactly" of the  form (\ref{Estar}).  
\begin{itemize}
\item \underbar{The case $n=1$.} 
When $n=1$, both nested sequences (\ref{twonested}) reduce to $\E_1$, for which the critical length is $\ell_1=\pi$. The condition $b-a<\pi$ is thus the necessary and sufficient condition for the existence of a one-dimensional W-space $\E_0^\star$ on $[a,b]$ contained in $\E_1$. We know that  there are  infinitely many such spaces, corresponding to inclusions of the form
$$
\E_0^\star=EC(\omega_0)\subset \E_1^\star=\E_1=EC(\omega_0, \omega_1).
$$ 
According to Theorem \ref{th:w0} and Remark  \ref{effectonbases} they all are obtained with 
\begin{equation}
\label{omega}
\omega_0:=\alpha_0\beta_0+\alpha_1\beta_1\hbox{ for some positive }\alpha_0, \alpha_1, \qquad \omega_1:=\left(\frac{\alpha_1\beta_1}{\omega_0}\right)', 
\end{equation}
where $(\beta_0, \beta_1)$ denotes the Lagrange basis of $\E_1$ relative to $(a,b)$,  that is
$$
\beta_0(x):=\frac{\sin(b-x)}{\sin(b-a)}, \quad \beta_1(x):=\frac{\sin(x-a)}{\sin(b-a)}, \quad x\in[a,b]. 
$$
\item \underbar{The left nested sequence in (\ref{twonested})}.
Assume that $n\geq 2$. The left nested sequence can be completed into a sequence of the form (\ref{Estar}) if and only we can find an $n$-dimensional W-space on $[a,b]$, $\F_{n-1}$  such that $\PP_{n-2}\subset \F_{n-1}\subset \E_n$. This is known  to be possible if and only if the two-dimensional space $D^{n-1}\E_n$ is an EC-space on $[a,b]$, see \cite{BO}. Since $D^{n-1}\E_n=\E_1$, this requirement is equivalent to $b-a<\ell_1=\pi$. Suppose the latter condition to be fulfilled. We then have infinitely many different choices for the subspace $\F_{n-1}$, namely 
\begin{equation}
\label{first}
\F_{n-1}=EC(\underbrace{\un, \un, \ldots, \un}_{ (n-2)\hbox{ times} }, \omega_0) \ \subset \ 
\E_{n}=EC(\underbrace{\un, \un, \ldots, \un}_{ (n-2)\hbox{ times} }, \omega_0,  \omega_1),
\end{equation}
where $\omega_0, \omega_1$ are defined in (\ref{omega}). The corresponding nested sequence (\ref{Estar}) is therefore given by 
$\E_i^\star=EC(w_0, \ldots, w_i)$, $i=0, \ldots, n$, with
\begin{equation}
\label{firstex}
w_i:= \un\hbox{ for }i=0, \ldots, n-2, \quad w_{n-1}:=\omega_0, \quad w_{n}:=\omega_1. 
\end{equation}
Since $b-a<\pi=\ell_1$, for each $n\geq 2$, the space $\E_{n}$ is an EC-space good for design on $[a,b]$.
Relative to $(a,b)$ we can thus compute the Bernstein basis in $\E_n$  via the integral recurrence relations formul\ae\ (\ref{IRR}) associated with the previous sequence of weight functions. 

\item \underbar{The right nested sequence in (\ref{twonested})}.
We now consider the right nested sequence in (\ref{twonested}). It can be completed into a sequence of the form (\ref{Estar}) if and only if we can find a one-dimensional EC-space $\E_0^\star \subset \E_1$, that is, once again if and only $b-a<\ell_1=\pi$. Supposing that $b-a<\pi$, through (\ref{omega}), the complete nested sequence $\E_0^\star \subset \E_1\subset \cdots \E_n$ provides us with another sequence $(w_0, \ldots, w_n)$ of weight functions on $[a,b]$, with $w_0:=\omega_0$, $w_1:=\omega_1$. As for $w_2, \ldots, w_n$ they are provided by the nested sequence $DL_1\E_2\subset \cdots\subset DL_1\E_{n-1}\subset DL_1\E_n$ (where the generalised derivative $L_1$ is defined by $\omega_0, \omega_1$), that is, after multiplication by $\omega_0\omega_1$, $\PP_0\subset\PP_1\subset \cdots\subset \PP_{n-2}$.  
We can therefore also write $\E_n=EC(w_0, \ldots, w_n)$ with
\begin{equation}
\label{secondex}
w_0:=\omega_0, \quad w_1:=\omega_1, \quad w_2=1/(\omega_0\omega_1), \quad w_i:= \un\hbox{ for }i=3, \ldots, n.
\end{equation}
\end{itemize}

We thus have at our disposal two obvious ways -- (\ref{firstex}) and (\ref{secondex}) -- to obtain the Bernstein basis relative to $(a,b)$ in $\E_n$ via integral recurrence relations (\ref{IRR}), and to use it to handle associated curves for design or any other purposes. Unfortunately, this is only valid under the assumption $b-a<\pi$. Most of the time,  trigonometric spaces are implicitly studied through (\ref{firstex}) for whatever purposes, for \eg \cite{ChenWang, MPS-R}. We can then avoid any reference to weight functions, starting the recurrence formul\ae\ (\ref{IRR}) from the Lagrange basis $(\beta_0, \beta_1)$ of $\E_1$. The same is generally done for the construction  for splines \cite{MP, wang, CosManPelSam, ManPelSam, ManPelSpe, ManReaSpe, IGA}.  

\medskip
\begin{figure}
\begin{center}
\includegraphics[height=0.22\textwidth]{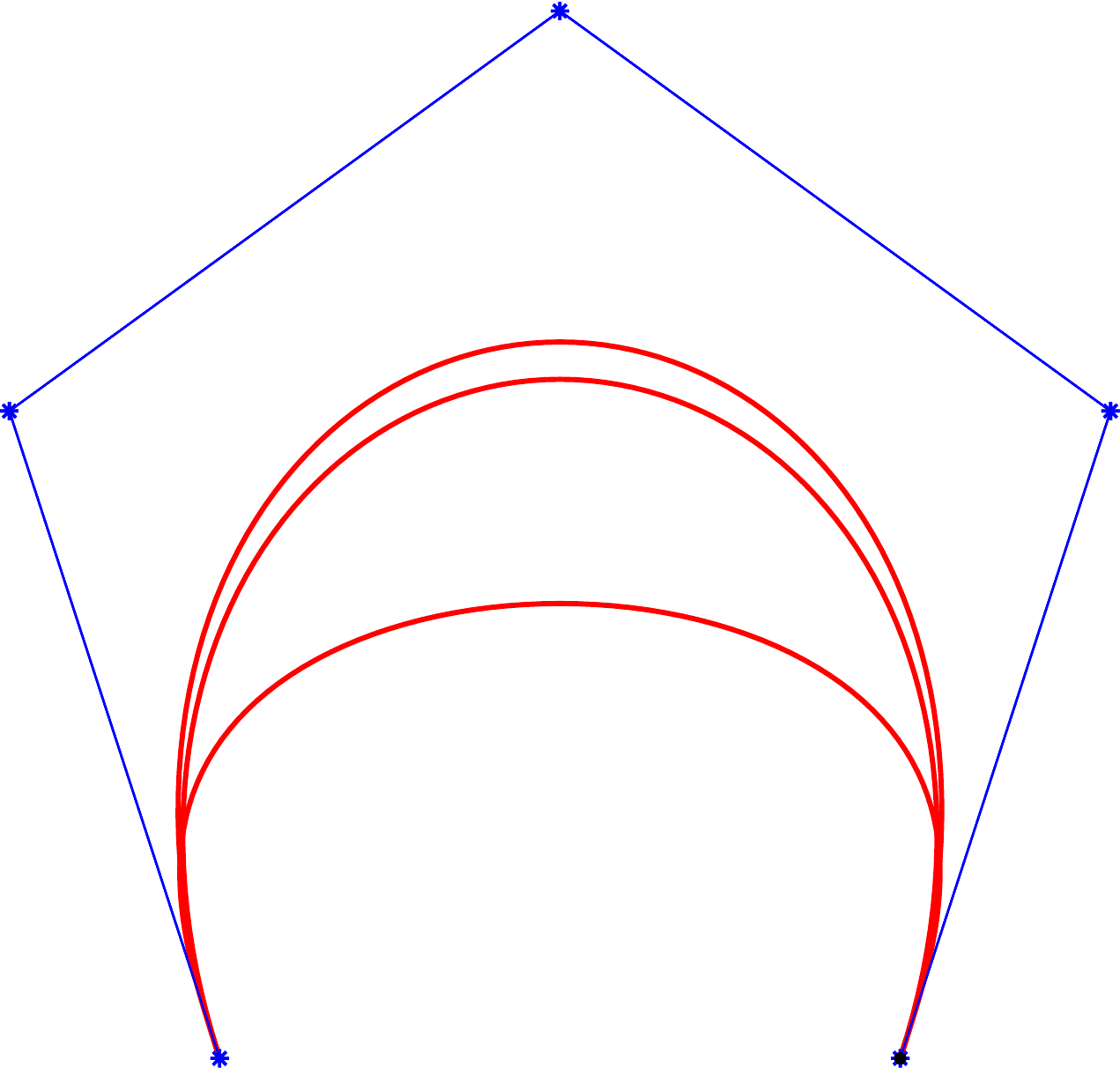} \hfill
\includegraphics[height=0.22\textwidth]{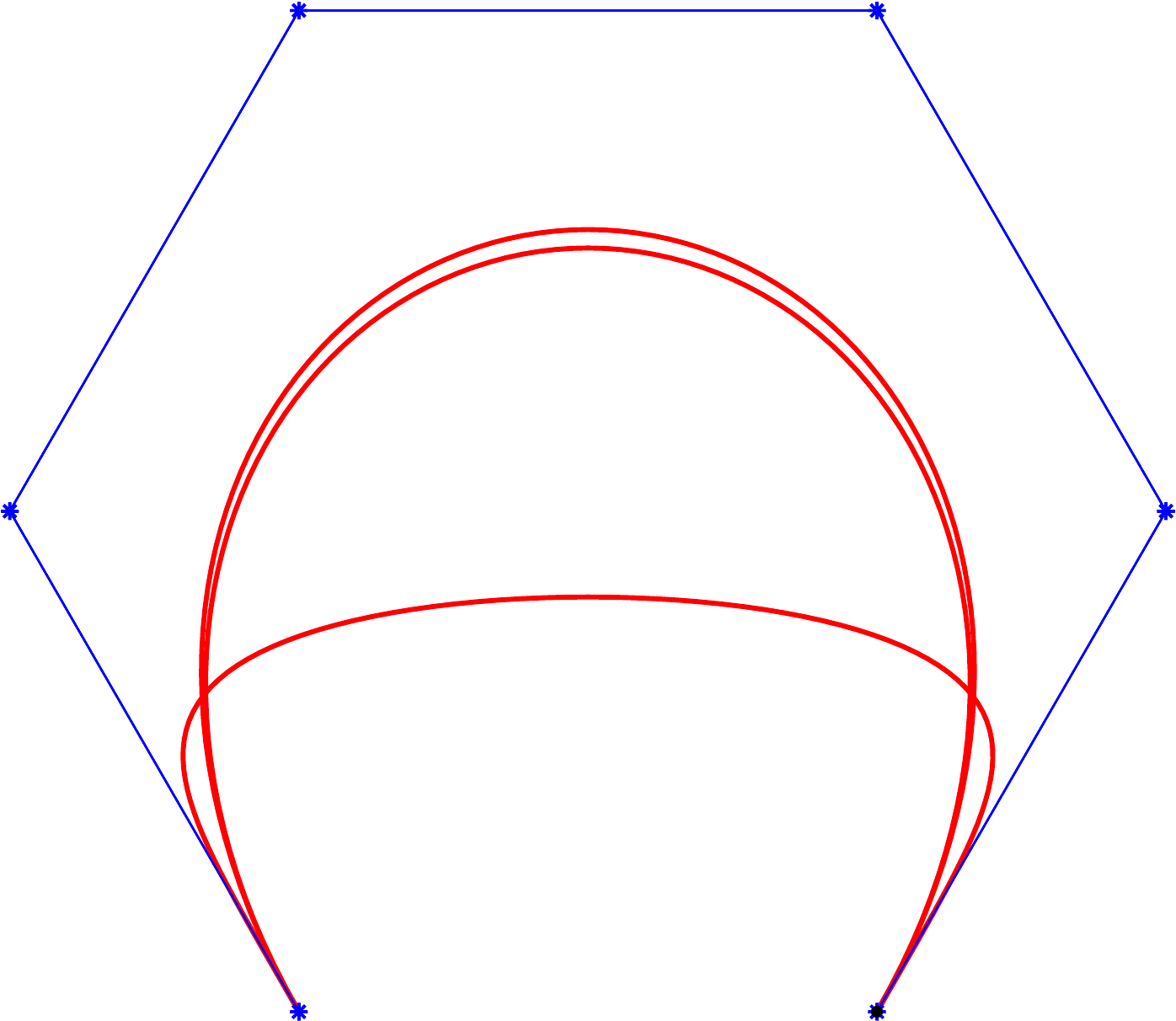} \hfill
\includegraphics[height=0.22\textwidth]{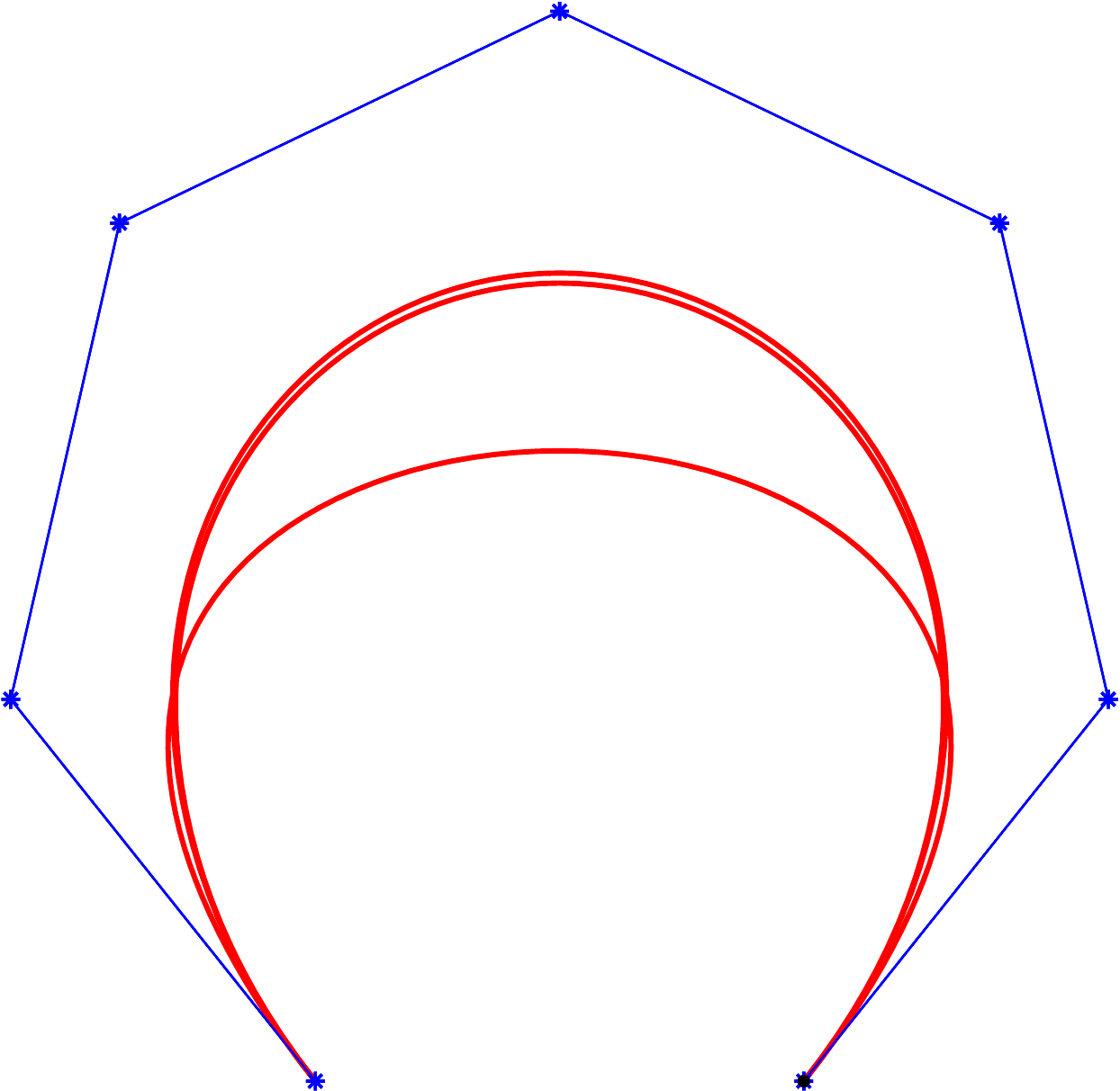} \hfill
\includegraphics[height=0.22\textwidth]{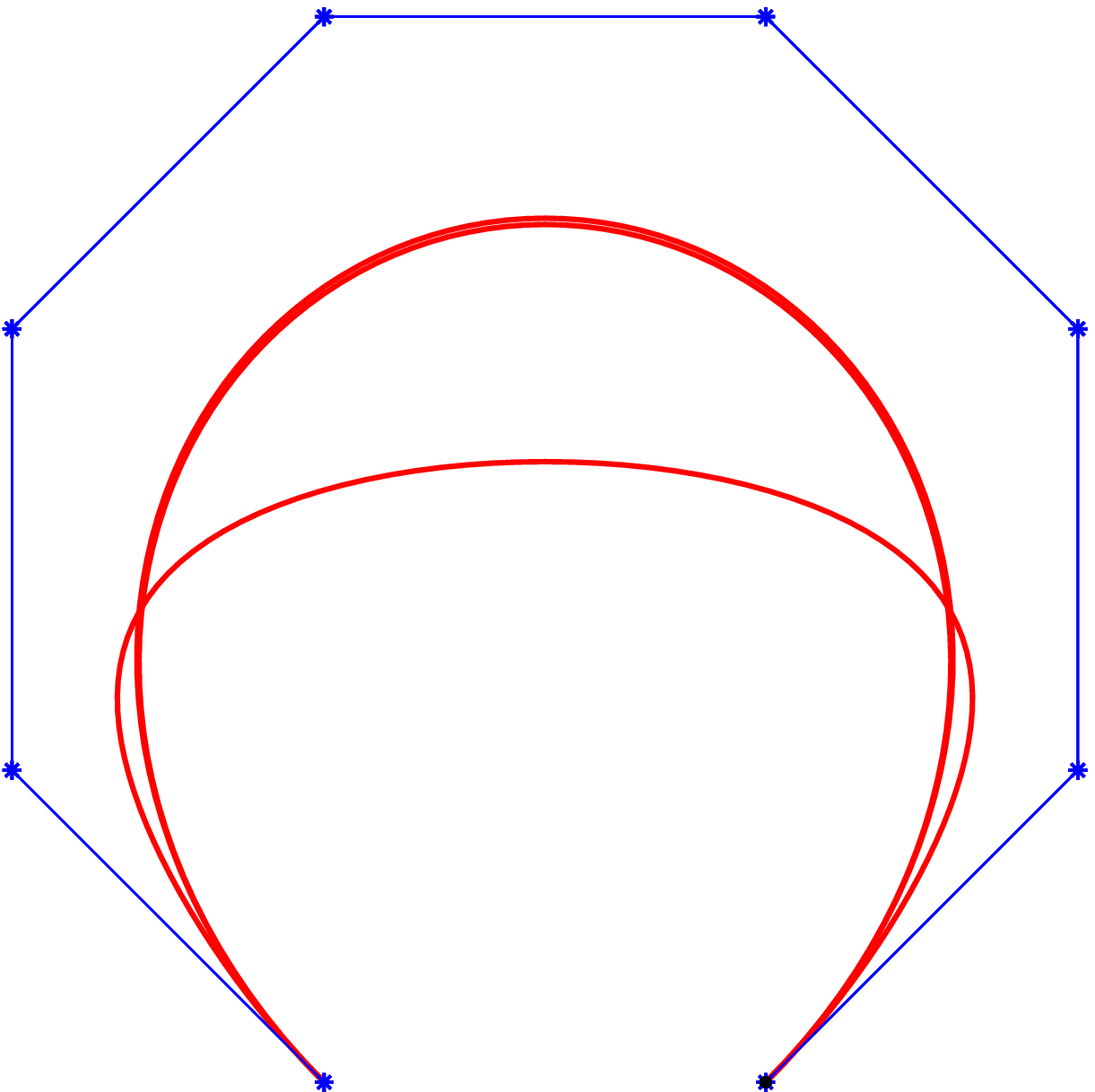} \hfill
\caption{For $n=4; 5; 6; 7$, Design in the space $\E_n$ spanned by $1, x, \ldots, x^{n-2}, \cos x, \sin x$ on $[0,h]$,  with, further and further from the control polygon, $h=0.1$; $h=3.14$; and $h=\ell_{n-1}^-$ (that is, $h=6.283$ for $n=4$; $h=8.986$ for $n=5,6$; $h=11.526$ for $n=7$).}
\label{fig:case4}
\end{center}
\end{figure}

As soon as $n\geq 3$, working under the assumption $b-a<\pi$  is actually a very limited use of trigonometric curves, if only with reference to (\ref{ellT}). To emphasise this point, observe that $\ell_n$ tends to infinity with $n$ due to (\ref{eq:T}). 
This limitation is illustrated in Figure 5  for increasing values of $n\geq 4$. A control polygon being given, we show a few corresponding curves in $\E_n$ on $[0, h]$ depending on $h$. When $h$ increases, they go from the polynomial curve of degree $n$ (visually obtained for $h=0.1$) to  the {\it critical curve} obtained for $h=\ell_{n-1}^-$  ($\ell_{n-1}$ being the critical length for design in the space  $\E_n$). In all cases we also show the curve obtained with  $h=3.14$, which is quite far from the critical curve, and can soon hardly be distinguished  from the polynomial curve as  $n$ increases. This naturally raises the following question: {\it for $n\geq 4$,  given that we simultaneously lose the remarkable simplicity of polynomials, is it really worthwhile replacing the degree $n$ polynomial space on $[a,b]$ by the trigonometric space $\E_n$ when requiring that $b-a<\pi$?} \footnote{Note that the case $n=3$ has been considered in many papers: it is well known that, when $h$ increases from $0^+$  to $2\pi^-$ the curve evolves from the cubic curve up to the segment joining the extreme control points (critical curve) and $h=\pi$  already yields a curve quite close to the polynomial one. Considering this case along with all pictures in Fig. 5, one can observe that the amplitude of the shape effects permitted by cycloidal spaces weakens in an ``oscillating" way as $n$ increases. This corresponds to the fact that the limit curves ignore the (two) central control point(s).} Under the weaker assumption $b-a<2\pi$, we similarly obtain the Bernstein basis relative to $(a,b)$ in the trigonometric space $\E_n$ on $[a,b]$ for any $n>3$ through formul\ae\ (\ref{IRR}) and ordinary integration, starting from the Bernstein basis in $\E_3$  (for which explicit expressions can be found in \cite{CAGD99}). Very soon as $n$ increases, we will be facing the same limitation. This can easily be guessed from the table giving the successive critical lengths  in \cite{CMP17}. 

This simple example clearly illustrates that only the knowledge of the critical length (for design) enables us to take full advantage of the parameter(s) attached to a given space of  functions. Not all kernels of differential operators have really appealing features, say for design. When investigating a new space, the most efficient way to learn if it produces remarkable shape effects and / or if it is worthwhile combining it with other spaces to build splines, is to consider its critical curves, that is, the curves obtained close to the critical length for design, which should therefore be determined beforehand, see \cite{BM12, criterion, JCAM18}. 
Unfortunately, we cannot expect this to be done experimentally from the visual analysis of either the curves or the Bernstein-type basis. Of course, if a parametric curve visibly contradicts the shape of its control polygon, or if one function of the expected Bernstein basis clearly takes negative values, we are certainly beyond the critical length for design. However, it is necessary to stress that we can be beyond it in spite of a visually satisfying behaviour. Moreover, the existence of a Bernstein basis / of a normalised totally positive basis implies that we are within the critical length for design only if, beforehand, we know that we are within the critical length (see Corollary \ref{cor:w0} and Theorem 4.1 of  \cite{CMP04}). 

These comments clearly emphasise how important it is to have at our disposal a reliable numerical method for the computation of critical lengths.   We would like to mention  that our procedure can also be useful on the theoretical side. It is well known that most results of the polynomial framework extend to  EC-spaces. Nonetheless, there are a few exceptions. To point out such exceptions, it can be necessary to observe what happens for kernels of linear differential operators with constant coefficients sufficiently close to their critical lengths. This will be illustrated in future work concerning dimension elevation.

\medskip

\noindent{\bf Acknowledgements:} The first two authors gratefully acknowledge support from
INdAM-GNCS Gruppo Nazionale per il Calcolo Scientifico.

\small {
}
\end{document}